\documentclass[draft]{amsart}
\usepackage{amsmath,amsfonts,amsrefs,latexsym}

\newtheorem{theorem}{Theorem}

\newtheorem{lemma}{Lemma}

\begin{document}
\title{The Vlasov-Maxwell-Boltzmann System in The Whole Space}
\author{Robert M. Strain}
\address{Department of Mathematics \\ Harvard University \\  
Cambridge \\ MA 02138 \\ USA} 
\email{strain at math.harvard.edu}
\subjclass[2000]{Primary: 76P05; 
Secondary: 82B40, 82C40, 82D05}

\begin{abstract}
The Vlasov-Maxwell-Boltzmann system is a fundamental model to describe the dynamics of dilute charged particles, where particles interact via collisions and through their self-consistent electromagnetic field.  We prove the existence of global in time classical  solutions to the Cauchy problem near Maxwellians.  
\end{abstract}

\maketitle

\thispagestyle{empty}

\section{The Vlasov-Maxwell-Boltzmann System}

The Vlasov-Maxwell-Boltzmann system is a very fundamental model to describe the 
dynamics of  dilute charged particles (e.g. electrons and ions):
\begin{equation}
{ \small
\begin{split}
\partial _tF_{+}+v\cdot \nabla _xF_{+}
+
\frac{e_{+}}{m_{+}}\left(E+\frac{v}{c}\times B\right)\cdot \nabla _vF_{+}
&=
Q^{++}(F_{+},F_{+})+Q^{+-}(F_{+},F_{-}),
 \\
\partial _tF_{-}+v\cdot \nabla _xF_{-}
-
\frac{e_{-}}{m_{-}}\left(E+\frac{v}{c}\times B\right)\cdot \nabla _vF_{-}
&=
Q^{--}(F_{-},F_{-})+Q^{-+}(F_{-},F_{+}).
\label{vlasov} 
\end{split}}
\end{equation}
The initial conditions are $F_{\pm }(0,x,v)=F_{0,\pm }(x,v)$. 
Here $F_{\pm }(t,x,v)\ge 0$ are number density
functions for ions (+) and electrons (-) respectively at time $t\ge 0$,
position $x=(x_1,x_2,x_3)\in \mathbb{R}^3$
and velocity $v=(v_1,v_2,v_3)\in {\mathbb R}^3$. 
The constants $e_{\pm }$ and $m_{\pm }$ are the magnitude of
their charges and masses, and $c$ is the speed of light.

The self-consistent electromagnetic field 
$[E(t,x),B(t,x)]$ in (\ref{vlasov}) is coupled with $F(t,x,v)$ through the celebrated
Maxwell system: 
\begin{equation*}
\begin{split}
\partial _tE-c\nabla_x \times B
&=
-4\pi \mathcal{J}=-4\pi \int_{{\mathbb R}^3}v\{e_{+}F_{+}-e_{-}F_{-}\}dv, 
\\
\partial _tB+c\nabla_x \times E &=0,
\end{split}
\end{equation*}
with constraints
\begin{equation*}
\begin{split}
\nabla_x \cdot B=0, 
~~~~
\nabla_x \cdot E=4\pi \rho =4\pi \int_{{\mathbb R}^3}\{e_{+}F_{+}-e_{-}F_{-}\}dv,
\end{split}
\end{equation*}
and initial conditions 
$
E(0,x)=E_0(x), B(0,x)=B_0(x).
$

Let $g_A(v),$ $g_B(v)$ be two number density functions for two types of particles $A$ and $B$ with masses $m_i$ and diameters $\sigma _i$ ($i\in\{A,B\}$).  In this article we consider the Boltzmann collision operator with hard-sphere interactions \cite{MR0258399}:
\begin{equation}
{\small
\begin{split}
Q^{AB}(g_A,g_B)
=
\frac{(\sigma _A+\sigma _B)^2}{4}\int_{{\mathbb R}^3\times S^2}|(u-v)\cdot \omega |
\{g_A(v^{\prime })g_B(u^{\prime })-g_A(v)g_B(u)\}dud\omega.
\label{hard} 
\end{split}}
\end{equation}
Here $\omega \in S^2$ and the post-collisional velocities are
\begin{equation}
v^{\prime }=v-\frac{2m_B}{m_A+m_B}[(v-u)\cdot \omega ]\omega ,\qquad
u^{\prime }=u+\frac{2m_A}{m_A+m_B}[(v-u)\cdot \omega ]\omega.  \label{prime}
\end{equation}
Then the pre-collisional velocities are $v,u$ and vice versa.

Recently global in time smooth solutions to  the Vlasov-Maxwell-Boltzmann system near Maxwellians were constructed by Guo \cite{MR2000470} in the spatially periodic case.   Then convergence  to Maxwellian for large times with any polynomial decay rate was shown for both the Vlasov-Maxwell-Boltzmann system  and the relativistic Landau-Maxwell system \cite{MR2100057}  by Strain and Guo \cite{strainGUOalmost}  in the periodic box.  Yet the corresponding questions for the Cauchy problem have remained open.

A simpler model can be formally obtained when the speed of light is sent to infinity. This is the
 Vlasov-Poisson-Boltzmann system, where $B(t,x)\equiv 0$ and $E(t,x)=\nabla_x \phi(t,x)$ in 
\eqref{vlasov}.  Here the self-consistent electric potential, $\phi$, satisfies a Poisson equation.   
There are several results for this model.

  Renormalized solutions of large amplitude were constructed by Lions \cite{MR1296258}, and this was  generalized to the case with boundary by Mischler \cite{MR1776840}.
The long time behavior of weak solutions with additional regularity assumptions are studied in Desvillettes and Dolbeault \cite{MR1104107}.  
The question of existence of 
renormalized solutions to the  Vlasov-Maxwell-Boltzmann system remains a major open problem in kinetic theory. 
 
Time decay of solutions to the linearized Vlasov-Poisson-Boltzmann system near Maxwellian was studied by Glassey and Strauss \cite{MR1669049,MR1696322}.  The existence of spatially periodic smooth solutions to the Vlasov-Poisson-Boltzmann system with near Maxwellian initial data was shown by Guo \cite{MR1908664}.    For some references on the Boltzmann equation with or without forces see \cite{MR1307620,MR0258399,MR1379589,MR1942465}.

Next we discuss global in time classical solutions to the Cauchy problem for the Vlasov-Poisson-Boltzmann system.  The case of near vacuum  data was solved by Guo \cite{MR1828983}  for some soft collision kernels.  
Duan, Yang and Zhu \cite{DYZ2005} have recently resolved the near vacuum data case for the hard sphere model and some cut-off potentials.
The near Maxwellian case was proven by Yang, Yu and Zhao  \cite{YYZ2004} with restrictions on either the size of the mean free path or the size of the background charge density.  After the completion of the results in this paper we learned that Yang and Yu  \cite{YZ2005} have removed those restrictions
 and further shown a time decay rate to Maxwellian of $O(t^{-1/2})$.  
 The method in \cite{YYZ2004,YZ2005} makes use of techniques from Liu and S.-H. Yu \cite{MR2044894} 
and Liu, Yang and S.-H. Yu \cite{MR2043729} 
 as well as the study in conservation laws. 
 But there seems to be serious obstacles to generalizing 
the method of \cite{YYZ2004,YZ2005}
to the full system \eqref{vlasov}  due to difficulties involved in controlling the electromagnetic field, in particular $B(t,x)$.

It is our purpose in this article to establish the existence of global in time classical solutions to the
Cauchy Problem for the Vlasov-Maxwell-Boltzmann system near global Maxwellians:
\begin{equation*}
\begin{split}
\mu _{+}(v)&=\frac{n_0}{e_{+}}
\left(\frac{m_{+}}{2\pi \kappa T_0}\right)^{3/2}
e^{-m_{+}|v|^2/2\kappa T_0},
\\
\mu _{-}(v)&=\frac{n_0}{e_{-}}\left(\frac{m_{-}}{2\pi \kappa T_0}\right)^{3/2}
e^{-m_{-}|v|^2/2\kappa T_0}. 
\end{split}
\end{equation*}
Accordingly, in the rest of this section we will reformulate the problem \eqref{vlasov} as a perturbation of the equilibria.

Since the presence of the physical constants do not create
essential mathematical difficulties, for notational simplicity,
we normalize all constants in the Vlasov-Maxwell-Boltzmann system to one. 
We normalize the Maxwellian
$$
\mu (v)\equiv \mu _{+}(v)=\mu _{-}(v)=(2\pi)^{-3/2}e^{-|v|^2/2}. 
$$
We further normalize the collision operator \eqref{hard} as
\begin{gather*}
Q^{AB}(g_A,g_B)
=
Q(g_A,g_B)
=
\int_{{\mathbb R}^3\times S^2}|(u-v)\cdot \omega |
\{g_A(v^{\prime })g_B(u^{\prime })-g_A(v)g_B(u)\}dud\omega.
\end{gather*}
And we define the standard perturbation $f_{\pm }(t,x,v)$ to $\mu $ as 
\begin{equation*}
F_{\pm }=\mu +\sqrt{\mu }f_{\pm }.  
\end{equation*}
Let $[\cdot ,\cdot ]$ denote a column vector so that 
$
F(t,x,v)=[F_{+},F_{-}]. 
$
We then study the normalized vector-valued Vlasov equation for the perturbation 
\[
f(t,x,v)=[f_{+}(t,x,v),f_{-}(t,x,v)],
\]
which now takes the form 
\begin{equation}
\{\partial _t+v\cdot \nabla _x+q(E+v\times B)\cdot \nabla _v\}f
-\{E\cdot v\}\sqrt{\mu }q_1+Lf=\frac{q}{2}\{E\cdot v\}f+\Gamma (f,f),  \label{vb}
\end{equation}
with $f(0,x,v)=f_0(x,v),$ $q_1=[1,-1]$
and 
$
q=
\text{diag}(1,-1)
$.

For any given $g=[g_+,g_-],$ the linearized collision operator in (\ref{vb}) is
$$
Lg=[L_+ g, L_- g],
$$ 
where
\begin{equation}
L_\pm g=
-2\mu^{-1/2}Q(\sqrt{\mu }g_\pm,\mu )
-\mu^{-1/2}Q(\mu ,\sqrt{\mu }\{g_\pm+g_\mp\}).  \label{L}
\end{equation}
We split $L$ in the standard way \cite{MR1379589}:
$
Lg = \nu(v) g -Kg.
$
The collision frequency is
\begin{equation}
\nu (v)\equiv  \int_{{\mathbb R}^3\times S^2} |v-u| \mu(u)du d\omega.
\label{nu}
\end{equation}
For $g=[g_+,g_-]$ and $h=[h_+,h_-]$, 
the nonlinear collision operator  is
$$
\Gamma (g,h)=[\Gamma_+ (g,h),\Gamma_- (g,h)],
$$
with 
\begin{equation*}
\Gamma_\pm (g,h)=
\mu^{-1/2} Q(\sqrt{\mu }g_\pm,\sqrt{\mu }h_\pm)
+\mu^{-1/2}Q(\sqrt{\mu }g_\pm,\sqrt{\mu }h_\mp).  \label{gamma}
\end{equation*}
Further, the coupled Maxwell system now takes the form 
\begin{eqnarray}
&&\partial _tE-\nabla_x \times B=-\mathcal{J}=-\int_{{\mathbb R}^3}v\sqrt{\mu }%
(f_{+}-f_{-})dv,\;\;\;\;\partial _tB+\nabla_x \times E=0,\;\;\;\;
\label{maxwell} \\
&&\nabla_x \cdot E=\rho =\int_{{\mathbb R}^3}\sqrt{\mu }(f_{+}-f_{-})dv,\;\;\;\;%
\nabla_x \cdot B=0,  \label{constraint}
\end{eqnarray}
with the same initial data $E(0,x)=E_0(x),\;B(0,x)=B_0(x).$

It is well known that the the linearized
collision operator $L$ is non-negative. For fixed $(t,x),$ the null
space of $L$ is given by (Lemma \ref{llinear}):
\begin{equation}
{\mathcal N}={\rm span}\{[\sqrt{\mu },0],\;[0,\sqrt{\mu }],\;[v_i\sqrt{\mu },v_i%
\sqrt{\mu }],\;[|v|^2\sqrt{\mu },|v|^2\sqrt{\mu }]\}~~~(1\le i\le 3).  \label{null}
\end{equation}
We define ${\bf P}$ as the orthogonal projection in $L^2({\mathbb R}^3_v)$ to the null
space ${\mathcal N}$. 
With $(t,x)$ fixed, we decompose any function 
$
g(t,x,v)=[g_+(t,x,v),g_-(t,x,v)]
$
as
\begin{equation*}
g(t,x,v)={\bf P}g(t,x,v)+({\bf I-P})g(t,x,v).  
\end{equation*}
Then ${\bf P}g$ is call the hydrodynamic part of $g$ and $({\bf I-P})g$ the
microscopic part.

\section{Notation and Main Results}

We shall use $\langle \cdot
,\cdot \rangle $ to denote the standard $L^2$ inner product in ${\mathbb R}^3_v$ for a pair of functions 
$g_1(v), g_2(v)\in L^2({\mathbb R}^3_{v};{\mathbb R}^2)$.  The corresponding norm is written 
$
|g_1|_2^2= \langle g_1, g_1\rangle.
$
We also define 
a weighted $L^2$ inner product in ${\mathbb R}^3_{v}$ as
\[
\langle g_1,g_2\rangle _\nu \equiv \langle \nu (v)g_1,g_2\rangle.
\]
And we use $|\cdot |_\nu $ for its corresponding $L^2$ norm. 
We note that the weight \eqref{nu} satisfies
$$
\frac{1}{c}(1+|v|)\le\nu(v)\le c(1+|v|).
$$
We also use $\|\cdot \| $ to denote $L^2$ norms in either ${\mathbb R}_x^3\times {\mathbb R}_v^3$ 
or ${\mathbb R}_x^3$.  In other words we use $\| \cdot \|$ to denote an $L^2$ norm in ${\mathbb R}_x^3\times {\mathbb R}_v^3$ if the function depends on $(x,v)$.  But if the function only depends upon 
$x\in {\mathbb R}_x^3$ then $\| \cdot \|$ will denote the norm of $L^2({\mathbb R}_x^3)$.  The  $L^2({\mathbb R}_x^3\times {\mathbb R}_v^3)$ inner product is written $(\cdot, \cdot)$.
Further, 
$\|g \|_\nu^2\equiv(\nu g, g)$.

Let multi-indices $\alpha$ and $\beta$ be 
\[
\alpha  =[\alpha _0,\alpha _1,\alpha _2,\alpha _3],\;\;
\beta =[\beta_1,\beta_2,\beta_3].
\]
We define a high order derivative as
\[
\partial _\beta ^\alpha \equiv 
\partial _t^{\alpha _0}\partial_{x_1}^{\alpha_1}\partial _{x_2}^{\alpha_2}\partial _{x_3}^{\alpha_3}
\partial _{v_1}^{\beta _1}\partial _{v_2}^{\beta _2}\partial_{v_3}^{\beta _3}. 
\]
If each component of $\alpha$ is not greater than that of $\bar{\alpha}$'s,
we write $\alpha \le \bar{\alpha};$ $\alpha <\bar{\alpha}$ means 
$\alpha \le \bar{\alpha}$ and $|\alpha |<|\bar{\alpha}|.$ We also denote 
$\begin{pmatrix}
\alpha \\ \bar{\alpha}
\end{pmatrix}$ by $C_{\bar{\alpha}}^\alpha .$ 

We next define an ``{\bf Instant Energy functional}''  for a solution to the Vlasov-Maxwell-Boltzmann system, $[f(t,x,v),E(t,x),B(t,x)]$, to be any functional $\mathcal{E}(t)$ which satisfies the following for an absolute constant $C>0$:
\begin{equation*}
\frac{1}{C}{\mathcal E}(t)\le 
\sum_{|\alpha |+|\beta |\le N}||\partial _\beta ^\alpha f(t)||^2
+\sum_{|\alpha |\le N}||[\partial ^\alpha
E(t),\partial ^\alpha B(t)]||^2
\le 
C
{\mathcal E}(t).
\end{equation*}
The temporal derivatives of $[f_0,E_0,B_0]$ are defined
through equations (\ref{vb}) and (\ref{maxwell}).  We will always use $C$ to denote a positive absolute  constant which may change from line to line.  

We additionally define the ``{\bf Dissipation rate}'' for a solution as 
\begin{equation*}
{\mathcal D}(t)\equiv 
\|E(t)\|^2+\sum_{0<|\alpha |\le N}||\partial^\alpha f(t)||^2_\nu+
\sum_{|\alpha |+|\beta |\le N}||\partial _\beta ^\alpha ({\bf I-P}) f(t)||^2_\nu.  
\end{equation*}
Notice that we include the electric field $E(t,x)$ in the dissipation and further we do not take any velocity derivatives of the hydrodynamic part, ${\bf P} f$.
Throughout this article we assume $N\ge 4$. Our main result is as follows.

\begin{theorem}
\label{main}Assume that $[f_0,E_0,B_0]$ satisfies the constraint (\ref{constraint})
initially. Let 
$$
F_{0,\pm }(x,v)=\mu +\sqrt{\mu }f_{0,\pm }(x,v)\ge 0.
$$ 
There exists an instant energy functional and $\epsilon_0>0$ such that if 
\[
{\mathcal E}{ }(0)\le \epsilon_0,
\]
then there is a unique global solution $[f(t,x,v),E(t,x),B(t,x)]$ to the
Vlasov--Maxwell-Boltzmann system (\ref{vb}) and (\ref{maxwell}) with (\ref
{constraint}) satisfying
$$
\mathcal{E}(t) +\int_0^t \mathcal{D}(s) ds\le \mathcal{E}(0).
$$
And moreover
$
F_{\pm }(t,x,v)=\mu +\sqrt{\mu }f_{\pm }(t,x,v)\ge 0.
$
\end{theorem}

In particular our solutions are $C^1$ for $N$ chosen large enough.
The proof of Theorem \ref{main} makes use of some of the  techniques developed in \cite{MR2000470,MR2095473}.  Yet several new  difficulties arise which are fundamental to the full space problem for the the Vlasov-Maxwell-Boltzmann system and we develop new techniques accordingly.

We expect that the techniques developed in this paper will be useful for working on other kinetic equations with force terms in the whole space.  In particular,  these techniques can be used to build solutions to the relativistic Landau-Maxwell system \cite{MR2100057} in the whole space.  
And our proof of Theorem \ref{main} supplies another construction of solutions to the Cauchy problem for the Vlasov-Poisson-Boltzmann system near Maxwellian.

As in \cite{MR2000470,MR2095473}, a key point in our construction is to show that the linearized collision operator \eqref{L} is effectively coercive for solutions of small amplitude to the Vlasov-Maxwell-Boltzmann system:

\begin{theorem}
\label{positive}Let $[f(t,x,v),E(t,x),B(t,x)]$ be a classical solution to (\ref{vb}) and (\ref{maxwell}) satisfying (\ref{constraint}). There exists $M_0>0$
such that if 
\begin{equation}
\sum_{|\alpha |\le N}\left\{||\partial ^\alpha f(t)||^2+||\partial
^\alpha E(t)||^2+||\partial ^\alpha B(t)||^2\right\}\le M_0,  \label{m0}
\end{equation}
then there are constants $\delta _0=\delta _0(M_0)>0$ and 
 $C_0>0$ so that
\[
\sum_{|\alpha |\le N}\left( L\partial ^\alpha f(t),\partial ^\alpha
f(t)\right) \ge \delta _0||({\bf I-P}) f(t)||_\nu^2+\delta _0\sum_{0<|\alpha |\le N}||\partial ^\alpha f(t)||_\nu
^2
-C_0\frac{dG(t)}{dt},
\]
where $G(t)= \int_{\mathbb{R}^3} \nabla _{x}\cdot b (a_++a_-)dx$ and $a_\pm, b$ are defined  below in  \eqref{p0}.
\end{theorem}

To prove Theorem \ref{positive}, we split the solution $f$ to the Vlasov-Maxwell-Boltzmann system (\ref{vb}) into it's hydrodynamic and microscopic parts as  
\begin{equation}
f={\bf P}f+({\bf I-P})f.
\label{microHYDROd}
\end{equation}
Separating into linear and
nonlinear parts, and using $L\{{\bf P}f\}=0$ (Lemma \ref{llinear}), we can express
the hydrodynamic part, ${\bf P}f$, through the microscopic part, $({\bf I-P})f,$ up to a
second order term: 
\begin{equation}
\{\partial _t+v\cdot \nabla _x\}{\bf P}f-\{E\cdot v\}\sqrt{\mu }q_1
=
l(({\bf I-P})f)+h(f),  \label{macro}
\end{equation}
where 
\begin{eqnarray}
l(({\bf I-P})f) &\equiv &-\{\partial _t+v\cdot \nabla _x+L\}({\bf I-P})f,  \label{l}
\\
h(f) &\equiv &-q(E+v\times B)\cdot \nabla _vf+\frac{q}{2}\{E\cdot v\}f+\Gamma (f,f).
\label{h}
\end{eqnarray}
Expanding ${\bf P}f=[{\bf P}_+f, {\bf P}_-f]$ as a linear combination of the basis in (\ref
{null}) yields
\begin{equation}
{\bf P}_\pm f
\equiv 
\left\{a_{\pm}(t,x)+\sum_{i=1}^3b_i(t,x)v_i+c(t,x)|v|^2\right\}\mu^{1/2}(v).
\label{p0}
\end{equation}
A system of ``macroscopic equations'' for  the coefficients
$a(t,x)=[a_{+}(t,x),a_{-}(t,x)]$, 
$b_i(t,x)$
 and
$c(t,x)$ can be derived from an expansion of the left side of (\ref{macro}) in the velocity variables, using
 \eqref{p0} to obtain  \cite{MR2000470}:
\begin{eqnarray}
\nabla _xc &=& l_c+h_c,  \label{c} \\
\partial _t c+\partial ^ib_i &=& l_i+h_i,  \label{bi} \\
\partial ^ib_j+\partial ^jb_i &=& l_{ij}+h_{ij},i\neq j,  
\label{bij} \\
\partial _t b_i+\partial ^ia_{\pm }\mp E_i
&=&
l_{bi}^{\pm }+h_{bi}^{\pm },
\label{ai} \\
\partial _t a_{\pm } &=& l_a^{\pm }+h_a^{\pm },  \label{adot}
\end{eqnarray}
Here $\partial ^j=\partial _{x_j}$.
The terms on the r.h.s. of the macroscopic equations \eqref{c}-\eqref{adot} are obtained 
by expanding the r.h.s. of \eqref{macro} with respect to the same velocity variables and comparing the coefficients on the two sides.
Here $l_c(t,x),l_i(t,x),$ $l_{ij}(t,x),$ $l_{bi}^{\pm }(t,x)$ and $l_a^{\pm }(t,x)$
are the coefficients of this expansion of the linear term 
\eqref{l}, 
and 
$h_c(t,x)$,
$h_i(t,x)$, 
$h_{ij}(t,x)$, 
$h_{bi}^{\pm }(t,x)$
and $h_a^{\pm }(t,x)$ are the coefficients of the same
expansion of the second order term \eqref{h}.  These terms are defined precisely in Section \ref{s:positivity}.

Theorem \ref{positive} is proven via energy estimates of these macroscopic equations \eqref{c}-\eqref{adot}, which are used to estimate derivatives of $a,b$ and $c$. 
In the periodic case \cite{MR2000470}, the Poincar{\'e} inequality was used to estimate $a,b$ and $c$ in terms of their derivatives and the conserved quantities.  
 Further, it was observed in \cite{MR2095473} that $a,b$ and $c$ themselves can not be estimated directly from \eqref{c}-\eqref{adot}.
 This is why we remove ${\bf P}f$ from the dissipation rate in the whole space.  
  Additionally, due to the velocity gradient in the equation \eqref{vb}, it is important to include the velocity derivatives in our norms.   On the other hand, in our global energy estimates (in Section \ref{s:global})
 $\partial_\beta {\bf P} f$ behaves almost the same as ${\bf P} f$.   We therefore take velocity derivatives  only on the microscopic part of the solution in the dissipation.  Our new weaker dissipation, $\mathcal{D}(t)$, allows the proof to work.    But at the same time this lack of control over the hydrodynamic part and it's velocity derivatives causes many of the  new difficulties.  
     
We use Sobolev type inequalities to bound the nonlinear terms by derivatives of $f$ in Lemma \ref{high}.  This is an important step for bounding the macroscopic equations by our weak dissipation in the proof of Theorem \ref{positive}.  Lemma \ref{high} is also used to include the electric field in the dissipation in Lemma \ref{field}.  
 
 Further, in the periodic case \cite{MR2000470} the Poincar{\'e} inequality was used to estimate the low order temporal derivatives of $b_i$.  Instead, we notice that there is a cancellation  of the most dangerous terms in the macroscopic equations which allows us to estimate the low order temporal derivatives of $b_i$ in the whole space directly in the proof of Theorem \ref{positive}.

Then we prove global existence via energy estimates of the Vlasov-Maxwell-Boltzman system in Theorem \ref{energyEST}.    New difficulties arise in these energy estimates for both the terms with  velocity derivatives and those without because we can not control the hydrodynamic part of a solution by the linear operator in the whole space.    We develop different procedures to prove the energy 
  estimates with the weaker dissipation.  In particular, we  estimate the nonlinear part of the collision operator  by the weaker dissipation in Lemma \ref{nonlin2}.  Further, including the electric field in the dissipation rate (Lemma \ref{field}) enables us to control other nonlinearities without velocity derivatives through Sobolev inequalities.  Additionally, we split the solution as \eqref{microHYDROd} in several places because this allows us to absorb velocity derivatives and exploit the exponential velocity decay  of the hydrodynamic part \eqref{p0} in our estimates.

This article is organized as follows: In section \ref{s:local} we recall some basic estimates for the Boltzmann collision operator.  In section \ref{s:positivity} we establish positivity of \eqref{L} for solutions to the Vlasov-Maxwell-Boltzmann system with small amplitude.  And in Section \ref{s:global} we prove  global existence (Theorem \ref{energyEST}).

\section{Local Solutions}\label{s:local}

 We will make use of the following basic estimates from 
\cite{MR1997264,MR1908664,MR2000470}:

\begin{lemma}[\cite{MR1997264,MR2000470}]
\label{llinear}$\langle Lg,h\rangle =\langle g,Lh\rangle ,\;\langle
Lg,g\rangle \ge 0.$ And $Lg=0$ if and only if 
$
g={\bf P}g,
$
where {\bf P} was defined in \eqref{p0}.
Furthermore there is a $\delta >0$ such that 
\begin{equation*}
\langle Lg,g\rangle \ge \delta |({\bf I-P})g|_\nu ^2.
\end{equation*}
\end{lemma}

From \eqref{L} it is well known that we can write $L=\nu-K$ where $\nu$ is the multiplication operator defined in \eqref{nu} and $K$ satisfies the following.

\begin{lemma}[\cite{MR1908664}]
\label{kcompact} If $\beta \neq 0,$ $\partial _\beta \nu (v)$ is uniformly bounded. For 
$g(v)\in H^{|\beta|} ({\mathbb R}^3;{\mathbb R}^2)$ and for any $\eta >0,$ there is $C_\eta >0$ such that 
\[
| \partial _\beta [Kg] |_2^2 \le \eta \sum_{|\beta^{\prime }|=|\beta |}|\partial _{\beta ^{\prime }}g|_2^2
+C_\eta |g|_2^2.
\]
\end{lemma}

Above and below $H^{|\beta|} ({\mathbb R}^3;{\mathbb R}^2)$ denotes the standard 
$L^2 ({\mathbb R}^3;{\mathbb R}^2)$ Sobolev space of order $|\beta|$.  Next up are estimates for the nonlinear part of the collision operator.

\begin{lemma}[\cite{MR1908664}]
\label{basic} Say $g(v),p(v), r(v)\in C^{\infty}_c ({\mathbb R}^3;{\mathbb R}^2)$. Then 
\begin{equation*}
\left| \langle \partial _\beta \Gamma (g_{,}p),r\rangle \right|   
\le 
C\sum_{ \beta _1+\beta _2\le \beta}
\left[ |\partial _{\beta _1}g|_\nu |\partial _{\beta_2}p|_2
+|\partial _{\beta _1}g|_2 |\partial _{\beta_2}p|_\nu\right] |r|_\nu.
\end{equation*}
Moreover,
\begin{equation*}
\left\| \int \Gamma (g_{,}p)rdv\right\| +\left\| \int \Gamma
(p_{,}g)rdv\right\| \le C\sup_{x,v}|\nu ^3r|\sup_x\left[ \int
|g(x,v)|^2dv\right] ^{1/2}||p||. 
\end{equation*}
\end{lemma}

Using these estimates, it is by now standard to prove existence and uniqueness of local-in time positive classical solutions to the Vlasov-Maxwell-Boltzmann system (\ref{vb}) and (\ref{maxwell}) with (\ref{constraint}) and to establish the continuity of the high order norms for a solution as in \cite[Theorem 4]{MR2000470}.  Therefore, we will not repeat the argument.

\section{Positivity of the Linearized Collision Operator}\label{s:positivity}

 In this section, we establish the positivity of \eqref{L} for any small solution 
$[f,E,B]$ to the  Vlasov-Maxwell-Boltzmann system (\ref{vb}) and (\ref{maxwell}) via a series of Lemma's.

First we simplify the Maxwell system from \eqref{maxwell} and \eqref{constraint}.  By (\ref{p0}) 
\begin{eqnarray*}
\int [v\sqrt{\mu },-v\sqrt{\mu }]\cdot {\bf P} fdv &=&0, \\
\int [\sqrt{\mu },-\sqrt{\mu }]\cdot fdv &=&a_{+}-a_{-}.
\end{eqnarray*}
We can therefore write the Maxwell system from \eqref{maxwell} and \eqref{constraint} as
\begin{equation}
\label{pf} 
\begin{split}
\partial _tE-\nabla_x \times B &=-\mathcal{J}
=
\int_{{\mathbb R}^3}[-v\sqrt{\mu },v\sqrt{\mu }]\cdot ({\bf I-P})fdv,
\\
\partial _tB+\nabla_x \times E&=0,
\\
\nabla_x \cdot E&=a_{+}-a_{-},
~~~\nabla_x \cdot B=0.
\end{split}
\end{equation} 
It is in this form that we will estimate the Maxwell system.

Next we discuss the derivation of the macroscopic equations \eqref{c}-\eqref{adot}.  We expand the r.h.s. of \eqref{macro} in the velocity variables using the projection  \eqref{p0} to obtain
\begin{gather*}
\left\{ 
v_i\partial ^ic|v|^2+\{\partial_t c+\partial
^ib_i\}v_i^2+\sum_{i<j}\{\partial ^ib_j+\partial ^jb_i\}v_iv_j+\{\partial_t
b_i+\partial ^ia_{\pm }\mp E_i\}v_i\right\} \sqrt{\mu }
\\
+\partial_t a_{\pm } \sqrt{\mu }, 
\end{gather*}
The above represents two equations.
Here $\partial ^j=\partial _{x_j}$.  For fixed ($t,x),$ this is an expansion of l.h.s. of \eqref{macro} with respect to the following basis ($1\le i<j\le 3)$: 
\begin{gather}
[v_i|v|^2\sqrt{\mu },v_i|v|^2\sqrt{\mu }],~
[v_i^2\sqrt{\mu },v_i^2\sqrt{\mu }], ~
[v_iv_j\sqrt{\mu },v_iv_j\sqrt{\mu }],
\nonumber
 \\
[v_i\sqrt{\mu },0],~
[0,v_i\sqrt{\mu }],~
[\sqrt{\mu },0],~
[0,\sqrt{\mu }].  
\label{base}
\end{gather}
By expanding the r.h.s. of \eqref{macro} with respect to the same velocity variables and comparing the coefficients of each basis element in \eqref{base} we obtain the macroscopic equations \eqref{c}-\eqref{adot}.   Given the full basis in \eqref{base}, for fixed $(t,x)$, $l_c(t,x),l_i(t,x),$ $l_{ij}(t,x),$ 
$l_{bi}^{\pm }(t,x)$ and $l_a^{\pm }(t,x)$, from \eqref{c}-\eqref{adot}, take the form 
\begin{equation}
\int_{{\mathbb R}^3}l(({\bf I-P})f)\cdot
\epsilon _n(v)dv.
\label{ldef}
\end{equation}
Here 
$l(({\bf I-P})f)$ is defined in \eqref{l} and
the $\{\epsilon _n(v)\}$,   which only depend on $v$, are linear combinations of the elements of \eqref{base}. Similarly 
$h_c$,
$h_i$,
$h_{ij}$,
$h_{bi}^{\pm}$ 
and $h_a^{\pm}$ from \eqref{c}-\eqref{adot} are also of the form 
\begin{equation}
\int_{{\mathbb R}^3} h(f)\cdot \epsilon _n(v)dv,
\label{hdef}
\end{equation}
where $h(f)$ is defined in \eqref{h}.  In the next two lemmas, we will estimate the right side of the macroscopic equations using \eqref{pf}, \eqref{ldef} and \eqref{hdef}.

\begin{lemma}
\label{linear}Let $\alpha =[\alpha _0,\alpha _1,\alpha _2,\alpha _3]$, $|\alpha |\le N-1$, then
for any $1\le i,j\le 3,$ 
\[
||\partial ^\alpha l_c||+||\partial ^\alpha
l_i||+||\partial ^\alpha l_{ij}||+||\partial ^\alpha l_{bi}^{\pm
}||+||\partial ^\alpha l_a^{\pm }||+||\partial ^\alpha \mathcal{J}||
\le 
C\sum_{|\bar{\alpha}|\le 1}||({\bf I-P})\partial^{\bar{\alpha}} \partial^\alpha f||.
\]
\end{lemma}

This Lemma is similar to \cite[Lemma 7]{MR2000470}, but we prove it for completeness.

\begin{proof}
The estimate for $\mathcal{J}$ follows from \eqref{pf}:
\begin{gather*}
|\partial ^\alpha \mathcal{J}| \le
C\int_{{\mathbb R}^3}|v|\sqrt{\mu }
\left|({\bf I-P})\partial ^\alpha f\right|dv 
\le 
C|({\bf I-P})\partial ^\alpha f|_2.
\end{gather*}
Square both sides and integrate over $x\in\mathbb{R}^3$ to get the estimate in Lemma \ref{linear} for $\mathcal{J}$.

For the other terms in Lemma \ref{linear}, it suffices to estimate the representation \eqref{ldef} with \eqref{l}.
Let $|\alpha |\le N-1.$
Since $\nu (v)\le C(1+|v|)$ and $K$ is bounded from 
$L^2({\mathbb R}^3;{\mathbb R}^2)$ to itself.   We have by (\ref{l}) 
\begin{gather*}
\left(\int \partial ^\alpha l(({\bf I-P})f)\cdot \epsilon _n(v)dv\right)^2
=
\left(\int (\{\partial _t+v\cdot \nabla _x+L\}({\bf I-P})\partial ^\alpha f)\cdot \epsilon _n(v)dv\right)^2, 
\\
\le 
C\int_{{\mathbb R}^3} |\epsilon _n(v)|dv 
 \int_{{\mathbb R}^3}|\epsilon_n(v)|(|({\bf I-P})\partial_t\partial ^\alpha f|^2+|v|^2|({\bf I-P})\nabla
_x\partial ^\alpha f|^2) dv
\\
+C
\int_{{\mathbb R}^3} |\epsilon _n(v)|dv 
 \int_{{\mathbb R}^3}|\epsilon_n(v)||(\nu -K)({\bf I-P})\partial ^\alpha f|^2dv
\\
\le 
C\{|({\bf I-P})\partial _t\partial ^\alpha f|^2_2+|({\bf I-P})\nabla _x\partial^\alpha f|_2^2
+|({\bf I-P})\partial ^\alpha f|_2^2\}.
\end{gather*}
Here we have used
$\partial^\alpha ({\bf I-P})f=({\bf I-P})\partial^\alpha f,$ and the exponential decay of 
$\epsilon _n(v)$.  The estimate in Lemma \ref{linear} follows by further integrating over $x\in\mathbb{R}^3$.
\end{proof}

Next, we estimate coefficients of the higher order term $h(f)$ on the right side of the macroscopic equations \eqref{c}-\eqref{adot}.

\begin{lemma}
\label{high}Let (\ref{m0}) hold for $M_0>0.$ Then for 
$1\le i\ne j\le 3,$ 
\[
\sum_{|\alpha |\le N}\{||\partial ^\alpha h_c||+||\partial ^\alpha
h_i||+||\partial ^\alpha h_{ij}||+||\partial ^\alpha h_{bi}^{\pm
}||+||\partial ^\alpha h_a^{\pm }||\}
\le 
C\sqrt{M_0}\sum_{0<|\overline{\alpha} |\le N}||\partial^{\overline{\alpha}} f||.
\]
\end{lemma}

Notice that the left side of this estimate contains zero'th order derivatives but the right side does not.  We have to tighten the estimate in this way because ${\bf P}f$ (with no derivatives) is not part of the dissipation rate.   We use  Sobolev's inequality on the 
  second order terms to remove the terms without derivatives.

\begin{proof} 
Notice that it suffices to estimate \eqref{hdef} with \eqref{h}.
For the first term of $h(f)$ in (\ref{h}), we integrate by parts to get 
\begin{gather*}
\int \partial ^\alpha \{q(E+v\times B)\cdot \nabla _vf)\}\cdot \epsilon_n(v)dv 
=
\int \sum_j\partial ^\alpha \{q(E+v\times B)_j \partial_{v_j}f)\}\cdot \epsilon_n(v)dv 
\\
=\sum C^{\alpha}_{\alpha _1}
\int \partial_{v_j}
\{(\partial ^{\alpha_1}E+v\times \partial ^{\alpha _1}B)_jq\partial ^{\alpha -\alpha _1}f\}\cdot
\epsilon _n(v)dv 
\\
=
-\sum C^{\alpha}_{\alpha _1}\int (\partial ^{\alpha _1}E+v\times \partial^{\alpha _1}B)_j
\{q\partial ^{\alpha -\alpha _1}f\}\cdot \partial_{v_j}\epsilon _n(v)dv.
\end{gather*}
Take the square of the above and further integrate over $x\in\mathbb{R}^3$, the result is
\begin{gather}
\le 
C\sum\int_{{\mathbb R}^3} \{|\partial ^{\alpha _1}E|^{2}+|\partial ^{\alpha _1}B|^{2}\}\left\{
\int |\partial ^{\alpha -\alpha _1}f|^2dv\right\}dx.
\label{mid}
\end{gather}
This follows from Cauchy-Schwartz and the exponential decay of $\nabla _v\epsilon _n(v)$. 

We will estimate \eqref{mid} in two steps.
First assume $|\alpha-\alpha_1|\le N/2$.  
In this case we use the embedding 
$W^{1,6}({\mathbb R}^3)\subset L^\infty ({\mathbb R}^3)$ combined with the $L^6({\mathbb R}^3)$ Sobolev inequality for gradients to obtain
\begin{equation*}
\begin{split}
\sup_{x\in\mathbb{R}^3} \int_{{\mathbb R}^3} |\partial ^{\alpha -\alpha _1} f(x,v)|^2dv 
&\le 
C\int_{{\mathbb R}^3} \|\partial ^{\alpha -\alpha _1} f\|^2_{W^{1,6}({\mathbb R}^3)}dv, 
\\
&\le 
C\int_{{\mathbb R}^3} \|\nabla_x\partial ^{\alpha -\alpha _1} f\|^2_{W^{1,2}({\mathbb R}^3)}dv.
\end{split}
\end{equation*}
Therefore, when
$
|\alpha-\alpha_1|\le N/2
$
and $N\ge 4$
we see that 
\begin{equation*}
\sqrt{\left| \eqref{mid} \right|}
\le
C
\{||\partial^{\alpha_1} E(t)||+||\partial^{\alpha_1} B(t)||\}
\sum_{0<|\bar{\alpha} |\le 2}||\partial^{\bar{\alpha}}  \partial ^{\alpha -\alpha _1} f(t)||.
\end{equation*}
The small amplitude assumption (\ref{m0}) completes the estimate in this case.

Alternatively if
$
|\alpha-\alpha_1|>N/2,
$
then $|\alpha _1|\le N/2$ and we use  the embedding $H^2({\mathbb R}^3)\subset L^\infty ({\mathbb R}^3)$ to obtain 
\[
\sup_{x\in {\mathbb R}^3}
\left(|\partial ^{\alpha _1}E|+|\partial ^{\alpha _1}B|\right)\le
C\sum_{|\bar{\alpha}|\le 2}\{||\partial^{\bar{\alpha}} \partial ^{\alpha _1} E(t)||
+||\partial^{\bar{\alpha}} \partial^{\alpha _1} B(t)||\}. 
\]
In this case
\begin{equation*}
\sqrt{\left| \eqref{mid} \right|}
\le
C
||\partial ^{\alpha-\alpha_1} f||\sum_{|{\bar{\alpha}}|\le 2}\{||\partial ^{\bar{\alpha}} \partial ^{\alpha _1} E(t)||
+||\partial ^{\bar{\alpha}} \partial ^{\alpha _1} B(t)||\}.
\end{equation*}
Again, the small amplitude assumption (\ref{m0}) completes the estimate.

For the second term of $h(f)$ in (\ref{h}) we have
\begin{gather*}
\int \frac{q}{2}\partial ^\alpha \{(E\cdot v)f\}\cdot \epsilon _n(v)dv 
=
\sum C^{\alpha}_{\alpha _1}\int \left\{\frac{q}{2}(\partial ^{\alpha _1}E\cdot v)\partial
^{\alpha -\alpha _1}f\right\}\cdot \epsilon _n(v)dv \\
\le C\sum |\partial ^{\alpha _1}E|\left\{ \int |\partial ^{\alpha -\alpha
_1}f|^2dv\right\} ^{1/2}.
\end{gather*}
This term is therefore easily treated by the last argument.

For the third term of $h(f)$, we have
\begin{gather*}
\left\| \int \partial ^\alpha \Gamma (f,f)\cdot \epsilon _n(v)dv\right\| 
\le 
\sum C^{\alpha}_{\alpha _1}\left\| \int \Gamma (\partial ^{\alpha_1}f,\partial ^{\alpha -\alpha _1}f)\cdot \epsilon _n(v)dv\right\|.
\end{gather*}
Without loss of generality assume $|\alpha_1|\le N/2$.  From Lemma \ref{basic}, the last line is
\begin{gather*}
\le C\sup_{x\in\mathbb{R}^3}\left\{ \int |\partial ^{\alpha _1}f(t,x,v)|^2dv\right\}^{1/2} ||\partial ^{\alpha-\alpha_1}f(t)||.
\end{gather*}
Again, we first use the embedding 
$W^{1,6}({\mathbb R}^3)\subset L^\infty ({\mathbb R}^3)$
and second the $L^6({\mathbb R}^3)$ Sobolev inequality to see that the above is
\begin{gather*}
\le 
C ||\partial ^{\alpha-\alpha_1}f(t)||\sum_{0<|\bar{\alpha} |\le 2}||\partial^{\bar{\alpha}}\partial ^{\alpha _1} f(t)||
\le C\sqrt{M_0}\sum_{0<|\bar{\alpha} |\le 2}||\partial^{\bar{\alpha}}\partial ^{\alpha _1} f(t)||.
\end{gather*}
The last inequality follows from (\ref{m0}). This completes the estimate of $h(f)$.
\end{proof}

Next we estimate the electric field $E(t,x)$ in terms of $f(t,x,v)$ through the macroscopic equation (\ref{ai}) using Lemma \ref{high}.

\begin{lemma}
\label{field}
Say $[f,E,B]$ is a classical solution to (\ref{vb}) and (\ref
{maxwell}) with (\ref{constraint}).  Let
the small amplitude assumption (\ref{m0}) be valid for some $M_0\le 1$. 
 $\exists C>0$ such that 
\[
\frac{1}{C}\sum_{|\alpha |\le N-1}||\partial ^\alpha E(t)||
\le 
\| ({\bf I-P})f\|+\sum_{0<|\alpha |\le N}||\partial ^\alpha f(t)||.
\]
\end{lemma}

Lemma \ref{high} is used to remove the hydrodynamic part, ${\bf P}f$, from this upper bound.

\begin{proof} We use the plus part of the macroscopic equation (\ref{ai}):
\[
-2\partial ^\alpha E_i
=
\partial ^\alpha l_{bi}^{+}+\partial ^\alpha h_{bi}^{+}
-
\partial ^\alpha \partial_t b_i
-
\partial ^\alpha \partial^i a_{+}. 
\]
By \eqref{p0}, 
\[
\|\partial ^\alpha \partial_tb_i\|
+
\|\partial ^\alpha \partial ^i a_{+}\|
\le
C\{\|{\bf P}\partial ^\alpha \partial_t f\|
+
\|{\bf P}\partial ^\alpha \partial^i f\|\}
\le 
C\{\|\partial ^\alpha \partial_t f||+\|\partial ^\alpha \partial^i f\|\}. 
\]
Since $M_0\le 1$ in assumption (\ref{m0}), applying Lemma \ref{linear} to $\partial ^\alpha l_{bi}^{+}$ and Lemma \ref{high} to $\partial ^\alpha h_{bi}^{+}$ we deduce Lemma \ref{field}.
\end{proof}

We now prove the positivity of \eqref{L} for a classical solution of small amplitude  to the full Vlasov-Maxwell-Boltzmann system, $[f(t,x,v),E(t,x),B(t,x)]$. The proof makes use of the estimates established in this section (Lemma \ref{linear}, \ref{high}, \ref{field}) in addition to an exact cancellation property of the macroscopic equations. \\

\noindent {\it Proof of Theorem} \ref{positive}. From Lemma \ref{llinear} we have
\[
\langle L\partial ^\alpha f,\partial ^\alpha f\rangle \ge \delta
|({\bf I-P})\partial ^\alpha f|_\nu ^2.
\]
It thus suffices to show that if (\ref{m0}) is valid for some small $M_0>0$,
then there  are constants $C_1, C_2>0$ such that 
\[
\sum_{0<|\alpha |\le N}\|{\bf P}\partial ^\alpha f(t)||_\nu^2 
\le 
C_1\sum_{|\alpha |\le N}||({\bf I-P})\partial^\alpha f(t)||_\nu ^2
+
C_2\frac{dG(t)}{dt}.
\]
Here $G(t)$ is defined in Theorem \ref{positive}.
By \eqref{p0}, we clearly have
\[
\frac{1}{C}||{\bf P}\partial ^\alpha f(t)||_\nu^2 
\le 
||\partial^\alpha [a_{+},a_{-}]||^2+||\partial ^\alpha b||^2+||\partial^\alpha c||^2.
\]
The rest of the proof is therefore devoted to establishing 
\begin{eqnarray}
\sum_{0<|\alpha |\le N}\{||\partial ^\alpha [a_{+},a_{-}]||^2
+||\partial ^\alpha b||^2+||\partial ^\alpha c||^2\}
&\le& C\sum_{|\alpha |\le N}
||({\bf I-P})\partial^\alpha f(t)||_\nu^2
\nonumber
\\
&&
+CM_0\sum_{0<|\alpha |\le N}||\partial ^\alpha f(t)||^2
\label{claim}
\\
&&
+
C_0\frac{dG(t)}{dt},
\nonumber
\end{eqnarray}
which implies Theorem \ref{positive} when $M_0$ is sufficiently small. 
This is because the last term on the right can be neglected for $%
M_0$ sufficiently small: 
\begin{eqnarray*}
||\partial ^\alpha f(t)||^2 
&=&
||{\bf P}\partial ^\alpha f(t)||^2+
||({\bf I-P})\partial^\alpha f(t)||^2 
\\
&\le &
\{\|\partial ^\alpha [a_{+},a_{-}]\|+||\partial^\alpha b||+||\partial ^\alpha c||\}^2
+
||({\bf I-P})\partial ^\alpha f(t)||^2.
\end{eqnarray*}
To prove (\ref{claim}), we estimate the macroscopic equations (\ref{c}) through (\ref{adot}).

\subsection*{The estimate for $b(t,x)$ with at least one spatial derivative}
We first estimate $\nabla_x \partial^\alpha b$ with 
$|\alpha |\le N-1$. 
We take $\partial ^j$ of (\ref{bi}) and (\ref{bij}) to get 
\begin{eqnarray*}
\Delta_x \partial ^\alpha b_i 
&=&\sum_j\partial ^{jj}\partial ^\alpha b_i 
=\left\{\sum_{j\ne i}\partial ^{jj}\partial ^\alpha b_i\right\}+\partial^{ii}\partial ^\alpha b_i 
\\
&=&
\sum_{j\neq i}\{-\partial ^{ij}\partial ^\alpha b_j+\partial ^j\partial
^\alpha l_{ij}+\partial ^j\partial ^\alpha h_{ij}\}
+
\{\partial ^i\partial^\alpha l_i+\partial ^i\partial ^\alpha h_i
-
\partial_t\partial ^i\partial^\alpha c\} 
\\
&=&
\sum_{j\neq i}\{\partial_t\partial ^i\partial ^\alpha c-\partial
^i\partial ^\alpha l_j-\partial ^i\partial ^\alpha h_j\}+\sum_{j\neq i}\{\partial ^j\partial ^\alpha l_{ij}+\partial ^j\partial
^\alpha h_{ij}\}
\\
&&
+\partial ^i\partial ^\alpha l_i+\partial ^i\partial ^\alpha
h_i 
-\partial_t\partial^i\partial ^\alpha c.
\end{eqnarray*}
Since 
$\sum_{j\neq i}\partial_t \partial^i\partial^\alpha c=2\partial_t\partial^i \partial^\alpha c$, we get
\begin{eqnarray*}
&=&
\partial_t\partial ^i\partial ^\alpha c 
+
\sum_{j\neq i}\{-\partial ^i\partial ^\alpha l_j-\partial ^i\partial
^\alpha h_j+\partial ^j\partial ^\alpha l_{ij}+\partial ^j\partial ^\alpha
h_{ij}\}+\partial ^i\partial ^\alpha l_i+\partial ^i\partial ^\alpha h_i 
\\
&=&
-\partial ^{ii}\partial ^\alpha b_i+\sum_{j\neq i}\{-\partial ^i\partial
^\alpha l_j-\partial ^i\partial ^\alpha h_j+\partial ^j\partial ^\alpha
l_{ij}+\partial ^j\partial ^\alpha h_{ij}\} 
+2\{\partial ^i\partial ^\alpha l_i+\partial ^i\partial ^\alpha h_i\}.
\end{eqnarray*}
Therefore, multiplying with $%
\partial ^\alpha b_i$ yields: 
\begin{equation}
||\nabla_x \partial ^\alpha b_i||\le C\sum \{||\partial ^\alpha
l_j||+||\partial ^\alpha h_j||+||\partial ^\alpha l_{ij}||+||\partial
^\alpha h_{ij}||+||\partial ^\alpha l_i||+||\partial ^\alpha h_i||\}.
\label{sb}
\end{equation}
This is bounded by the first two terms on the r.h.s. of (\ref{claim}) by Lemma \ref{linear} and 
Lemma \ref{high}. We treat the pure temporal derivatives of $b(t,x)$ at the
end the proof.

\subsection*{The estimate for $c(t,x)$}  From (\ref{c}) and (\ref{bi}), for $|\alpha |\le N-1$, we have 
\begin{eqnarray*}
||\nabla_x \partial ^\alpha c|| &\le &||\partial ^\alpha l_c||+||\partial
^\alpha h_c||,
\\
||\partial_t \partial ^\alpha c|| 
&\le &
||\partial ^i\partial ^\alpha b_i||+||\partial ^\alpha l_i||+||\partial ^\alpha h_i||.
\end{eqnarray*}
By (\ref{sb}), Lemma \ref{linear} and \ref{high}  both 
$||\partial_t \partial ^\alpha c||^2$ and 
$||\nabla_x \partial ^\alpha c||^2$ 
are bounded by the first two terms on the r.h.s. of (\ref{claim}).

\subsection*{The estimate for $[a_{+}(t,x),a_{-}(t,x)]$} 
By (\ref{adot}), for $|\alpha |\le N-1$, we get
\[
||\partial _t\partial ^\alpha [a_{+},a_{-}]||\le ||\partial ^\alpha l_a^{\pm
}||+||\partial ^\alpha h_a^{\pm }||.
\]
By Lemma \ref{linear} and \ref{high}, 
$||\partial _t\partial ^\alpha [a_{+},a_{-}]||$ is thus bounded by the first two terms on the r.h.s. of (\ref{claim}).

We now turn to purely spatial derivatives of $a_{\pm}(t,x)$.  Let $|\alpha |\le N-1$
and 
\[
\alpha =[0,\alpha _1,\alpha _2,\alpha _3]\neq 0. 
\]
By taking $\partial ^i$ of (\ref{ai}) and summing over $i$, we get 
\begin{equation}
-\Delta_x \partial ^\alpha a_{\pm }\pm \nabla_x \cdot \partial ^\alpha E
=
\nabla_x\cdot \partial_t \partial^\alpha b
-
\sum_i\partial ^i\partial ^\alpha\{l_{bi}^{\pm }+h_{bi}^{\pm }\}.  \label{div0}
\end{equation}
From the Maxwell system in (\ref{pf}) we have
\[
\nabla_x \cdot \partial ^\alpha E=\partial ^\alpha a_{+}-\partial^\alpha a_{-}. 
\]
We plug this into (\ref{div0}), then multiply (\ref{div0}) with $\partial ^\alpha a_{\pm }$ respectively and integrate over $\mathbb{R}^3_x$. 
Adding the two resulting equations yields
\begin{gather*}
||\nabla_x \partial ^\alpha a_{+}||^2
+
||\nabla_x \partial ^\alpha a_{-}||^2
+
||\partial^\alpha a_{+}-\partial ^\alpha a_{-}||^2 
\\
= 
\int \sum_{i, \pm} \partial^i  \partial^\alpha a_{\pm}  \left\{- \partial_t   \partial^\alpha b_i
+
\partial^\alpha l_{bi}^{\pm }+  \partial^\alpha h_{bi}^{\pm }
\right\} dx.
\end{gather*}
We therefore conclude that 
\[
||\nabla_x \partial ^\alpha a_{+}||+||\nabla_x \partial ^\alpha a_{-}||
\le
||\partial_t \partial ^\alpha b||
+
\sum_{i, \pm}||\partial ^\alpha \{l_{bi}^{\pm }+h_{bi}^{\pm }\}||. 
\]
Since $\alpha$ is purely spatial, this is bounded by the first two terms on the right side of (\ref
{claim}) by (\ref{sb}),  Lemma \ref{linear} and \ref{high}.

Now we consider the case $\alpha=0$. The same procedure yields 
\begin{equation*}
\frac{1}{2}||\nabla_x  a_{+}||^2
+
\frac{1}{2}||\nabla_x  a_{-}||^2
\leq 
\int \nabla _{x}\cdot \partial_t b(a_++a_-)+C\sum_{i, \pm}
||\{l_{bi}^{\pm }+h_{bi}^{\pm }\}||^{2}. 
\end{equation*}
We redistribute the derivative in $t$, as in \cite{MR2095473}, to estimate the first
term  by
\begin{eqnarray*}
\int \nabla _{x}\cdot \partial_t b(a_++a_-)dx 
&=&
\frac{d}{dt}\int \nabla _{x}\cdot b (a_++a_-)dx
-
\int \nabla _{x}\cdot b \partial_t (a_++a_-)  \label{a3} 
\\
&\leq &\frac{d}{dt}\int \nabla _{x}\cdot b(a_++a_-)dx
+C||\nabla_x b||^{2}+C\sum_\pm ||\partial_t a_\pm||^{2}.  \nonumber
\end{eqnarray*}
The first term on the right is $dG/dt$ and the last two terms have already been estimated by the other terms on the right side of \eqref{claim}.

\subsection*{The estimate for $b(t,x)$ with purely temporal derivatives}
Lastly, we consider $\partial_t \partial^\alpha b_i(t,x)$ with $\alpha =[\alpha _0,0,0,0]$ and $\alpha_0\le N-1$.   In \cite{MR2000470}, Guo used a procedure which involved the Poincar{\'e} inequality and Maxwell's equations \eqref{pf} to estimate $\partial^\alpha \partial_t b_i$ from the macroscopic equation \eqref{ai}, which includes the electric field.  But the Poincar{\'e} inequality does not hold in the whole space case to estimate the low order pure temporal derivatives of $b$.  

On the other hand, there are two equations for $\partial_t b_i$ in \eqref{ai}.  We  can add them together to cancel the most dangerous term, $E_i$, and obtain
$$
2\partial _t b_i+\partial ^ia_{+}+\partial ^ia_{-}
=
l_{bi}^{+}+h_{bi}^{+}+l_{bi}^{-}+h_{bi}^{-}.
$$
Without this exact cancellation, our estimates for $E_i$ are not sufficient to conclude the required estimates for 
$\partial_t \partial^\alpha b_i(t,x)$.
From this equation we can estimate the temporal derivatives of $b$ directly:
\begin{equation*}
||\partial_t \partial^\alpha b_i(t,x)||^2
\leq 
C\sum_\pm
\left(||\partial^i \partial ^\alpha a_{\pm}||^2
+
||\partial ^\alpha \{l_{bi}^{\pm}+h_{bi}^{\pm}\}||^2\right). 
\end{equation*}
If $\alpha_0>0$ then the terms on the right have already been estimated above by the first two terms terms on the right side of \eqref{claim}.

Alternatively if $\alpha_0=0$, we use the previous estimate for $\partial^i a_\pm$ to obtain
\begin{equation*}
\begin{split}
||\partial_t  b_i(t,x)||^2
\leq &
C\frac{d}{dt}\int \nabla _{x}\cdot b (a_++a_-)dx
+
C||\nabla_x b||^{2}
\\
&
+
C\sum_\pm \left( ||\partial_t a_\pm||^{2}+
|| \{l_{bi}^{\pm}+h_{bi}^{\pm}\}||^{2} \right). 
\end{split}
\end{equation*}
From Lemma \ref{linear} and \ref{high}, we finish Theorem \ref{positive} with $G=\int_{\mathbb{R}^3} \nabla _{x}\cdot b (a_++a_-)dx$.
\qed

\section{Global Existence}\label{s:global}

The main goal of this section is to establish global existence of classical solutions to the Vlasov-Maxwell-Boltzmann system, \eqref{vb} and \eqref{maxwell} with \eqref{constraint}.  Our proof uses Theorem \ref{positive} and energy estimates (Theorem \ref{energyEST}) developed in this section.

We have the following global existence theorem:

\begin{theorem}\label{energyEST}  Let $[f(t,x,v),E(t,x),B(t,x)]$ be a classical solution to 
to the Vlasov--Maxwell-Boltzmann system (\ref{vb}) and (\ref{maxwell}) with (\ref
{constraint}) satisfying \eqref{m0}.  Then there exists an instant energy functional satisfying
$$
\frac{d}{dt} \mathcal{E}(t)+\mathcal{D}(t)\le \sqrt{{\mathcal E}(t)}\mathcal{D}(t).
$$
With \eqref{m0}, this yields global existence from a standard continuity argument.  
\end{theorem}

First we prove the following nonlinear estimate:

\begin{lemma}\label{nonlin2}  Let $|\alpha|+|\beta|\le N$.  Then there is an instant energy functional such that
$$
\left| \left( \partial_\beta^\alpha \Gamma (f, f), \partial_\beta^\alpha ({\bf I-P})f\right) \right|
\le
\mathcal{E}^{1/2}(t)\mathcal{D}(t).
$$
\end{lemma}

The main improvement of this estimate over Lemma \ref{basic}  is the presence of the dissipation, $\mathcal{D}(t)$.
Clearly $\| E\|^2$ is not needed in the dissipation in Lemma \ref{nonlin2}.

\begin{proof}  First assume $\beta=0$.  We have
$$
\partial^\alpha \Gamma (f, f)=\sum_{ \alpha_1\le \alpha} C^\alpha_{\alpha_1}
\Gamma (\partial^{\alpha-\alpha_1}f ,\partial^{\alpha_1}f ).
$$
Furthermore 
(see for instance \cite[p.59-60]{MR1379589}):
\begin{equation}
\langle  \Gamma (g,p ),  r\rangle=\langle  \Gamma (g,p ), ({\bf I-P}) r\rangle.
\label{invariantP}
\end{equation}
Plugging these last two observations into Lemma \ref{basic} yields
\begin{equation*}
\left|  \langle \partial^\alpha \Gamma (f, f), \partial^\alpha f\rangle \right|   
\le 
C\sum_{ \alpha_1\le \alpha}
\left( |\partial^{\alpha-\alpha_1}f|_\nu |\partial^{\alpha_1}f|_2
+|\partial^{\alpha-\alpha_1}f|_2 |\partial^{\alpha_1}f|_\nu\right)
|({\bf I- P})\partial^{\alpha}f|_\nu.
\end{equation*}
Further integrate over $\mathbb{R}^3_x$ to obtain
\begin{gather*}
\left|  \left( \partial^\alpha \Gamma (f, f), \partial^\alpha ({\bf I-P}) f\right) \right|   
\\
\le 
C\sum_{ \alpha_1\le \alpha}
\int_{\mathbb{R}^3}\left( |\partial^{\alpha-\alpha_1}f|_\nu |\partial^{\alpha_1}f|_2
+|\partial^{\alpha-\alpha_1}f|_2 |\partial^{\alpha_1}f|_\nu\right)
|({\bf I-P})\partial^{\alpha}f|_\nu  dx.
\end{gather*}
By Cauchy-Schwartz and 
 symmetry 
the above is
\begin{gather*}
\le 
C\|({\bf I-P})\partial^{\alpha}f\|_\nu
\sum_{ \alpha_1\le \alpha}
\left\{\int_{\mathbb{R}^3}\left( |\partial^{\alpha-\alpha_1}f|_\nu |\partial^{\alpha_1}f|_2\right)^2
 dx\right\}^{1/2}.
\end{gather*}
Our goal is to show that 
$$
\left\{\int_{\mathbb{R}^3}\left( |\partial^{\alpha-\alpha_1}f|_\nu |\partial^{\alpha_1}f|_2\right)^2
 dx\right\}^{1/2}
\le \mathcal{E}^{1/2}(t)\mathcal{D}^{1/2}(t).
$$
This will imply Lemma \ref{nonlin2} in the case $\beta=0$.     We will use the embedding 
$W^{1,6}(\mathbb{R}^3)\subset L^{\infty}(\mathbb{R}^3)$
followed by the $L^{6}(\mathbb{R}^3)$ Sobolev inequality to one of the terms inside the $x$ integration.  
But, because the dissipation does not contain ${\bf P} f$,  it is not enough to just take the supremum in $x$ over the term with the smaller number of total derivatives. 
Instead, if  $|\alpha-\alpha_1|\le N/2$ then we always take the sup on the $\nu$ norm:
\begin{gather*}
\left\{\int_{\mathbb{R}^3}\left( |\partial^{\alpha-\alpha_1}f|_\nu |\partial^{\alpha_1}f|_2\right)^2
 dx\right\}^{1/2}
\le
\|\partial^{\alpha_1}f\|\sup_{x\in\mathbb{R}^3}  |\partial^{\alpha-\alpha_1}f(x)|_\nu
\\
\le
C\|\partial^{\alpha_1}f\| 
\left\{\int_{\mathbb{R}^3}\nu(v)\|\partial^{\alpha-\alpha_1}f(v)\|^2_{W^{1,6}(\mathbb{R}^3_x)}dv\right\}^{1/2}
\\
\le
C\|\partial^{\alpha_1}f\| \sum_{|\bar{\alpha}|\le 1} \left\{\int_{\mathbb{R}^3\times \mathbb{R}^3}
\nu(v)|\partial^{\alpha-\alpha_1}\partial^{\bar{\alpha}}\nabla_x f|^2dx dv\right\}^{1/2}
\le \mathcal{E}^{1/2}(t)\mathcal{D}^{1/2}(t).
\end{gather*}
Alternatively, if $|\alpha-\alpha_1|> N/2$ then $|\alpha_1|\le N/2$ and 
$\| \partial^{\alpha-\alpha_1}f\|_\nu$ is part of the dissipation so that 
it is ok to
take the supremum on the other norm to obtain
\begin{gather*}
\left\{\int_{\mathbb{R}^3}\left( |\partial^{\alpha-\alpha_1}f|_\nu |\partial^{\alpha_1}f|_2\right)^2
 dx\right\}^{1/2}
\le
\|\partial^{\alpha-\alpha_1}f\|_\nu \sup_{x\in\mathbb{R}^3}  |\partial^{\alpha_1}f(x)|_2
\\
\le
C\|\partial^{\alpha-\alpha_1}f\|_\nu
 \sum_{|\bar{\alpha}|\le 1} \left\{\int_{\mathbb{R}^3\times \mathbb{R}^3}
|\partial^{\alpha_1}\partial^{\bar{\alpha}}\nabla_x f|^2dx dv\right\}^{1/2}
\le \mathcal{D}^{1/2}(t)\mathcal{E}^{1/2}(t).
\end{gather*}
This completes the proof of Lemma \ref{nonlin2} when $\beta=0$.

If $|\beta|>0$ then
we expand $f$ as in \eqref{microHYDROd} to obtain
$$
\Gamma (f, f)=
\Gamma ({\bf P}f, {\bf P}f)
+
\Gamma ({\bf P}f, ({\bf I-P})f)
+
\Gamma (({\bf I-P})f, {\bf P}f)
+
\Gamma (({\bf I-P})f, ({\bf I-P})f).
$$
We will estimate each of these terms using similar Sobolev embedding arguments to above.
However, we use this splitting because $\partial_\beta {\bf P}f$ is not part of the dissipation.  The exponential velocity decay of ${\bf P}f$ is used to control the hard sphere weight \eqref{nu}.  Alternatively 
$({\bf I-P})f$ is part of the dissipation.
We only estimate the first two terms below.  The other two terms can be treated similarly.

Lemma \ref{basic}  and the definition \eqref{p0} yield 
\begin{gather*}
\left|  \langle \partial_\beta^\alpha \Gamma ({\bf P}f, {\bf P}f), \partial_\beta^\alpha ({\bf I-P})f\rangle \right|   
\\
\le 
C\sum_{ \alpha_1\le \alpha}
\left[ |\partial^{\alpha-\alpha_1} {\bf P}f|_2 |\partial^{\alpha_1} {\bf P}f|_2
+|\partial^{\alpha-\alpha_1} {\bf P}f|_2 |\partial^{\alpha_1}{\bf P}f|_2\right] 
|\partial_\beta^\alpha ({\bf I-P})f|_\nu.
\end{gather*}
Further integrate the above over $\mathbb{R}^3_x$, use Cauchy-Schwartz and symmetry to obtain
\begin{gather*}
\left|  \left( 
\partial_\beta^\alpha \Gamma ({\bf P}f, {\bf P}f), \partial_\beta^\alpha ({\bf I-P})f
\right) \right|   
\\
\le 
C\|\partial_\beta^\alpha ({\bf I-P})f\|_\nu
\sum_{\alpha_1\le \alpha}
\left\{\int_{\mathbb{R}^3}
|\partial^{\alpha-\alpha_1} {\bf P}f|_2^2 |\partial^{\alpha_1}{\bf P}f|_2^2
 dx\right\}^{1/2}.
\end{gather*}
Without loss of generality assume $|\alpha_1|\le N/2$.  Then by the same embedding
\begin{gather*}
\int_{\mathbb{R}^3}
|\partial^{\alpha-\alpha_1} {\bf P}f|_2^2 |\partial^{\alpha_1}{\bf P}f|_2^2
 dx
 \le
\left(\sup_{x\in\mathbb{R}^3}
|\partial^{\alpha_1}f(x)|_2^2\right)
\|\partial^{\alpha-\alpha_1} {\bf P}f\|^2 
\\
 \le
C\left(\sum_{|\bar{\alpha}|\le 1}
\|\partial^{\alpha_1}\partial^{\bar{\alpha}}\nabla_x f(x)\|^2\right)
\|\partial^{\alpha-\alpha_1} {\bf P}f\|^2 
\le \mathcal{D}(t)\mathcal{E}(t).
\end{gather*}
This yields the estimate in Lemma \ref{nonlin2} for the $\partial_\beta^\alpha \Gamma ({\bf P}f, {\bf P}f)$ term.

Finally, we estimate the $\partial_\beta^\alpha \Gamma ({\bf P}f, ({\bf I- P})f)$ term.  In this case
Lemma \ref{basic} yields
\begin{gather*}
\left|  \left( \partial_\beta^\alpha \Gamma ({\bf P}f,  ({\bf I- P}) f), \partial_\beta^\alpha ({\bf I-P})f\right) \right|   
\\
\le 
C\|\partial_\beta^\alpha ({\bf I-P}) f\|_\nu\sum
\left\{\int_{\mathbb{R}^3}\left(
|\partial _{\beta _1}^{\alpha-\alpha_1}{\bf P}f|_\nu |\partial _{\beta_2}^{\alpha_1}({\bf I-P})f|_2
\right)^2
 dx\right\}^{1/2}
 \\
 +
 C\|\partial_\beta^\alpha ({\bf I-P}) f\|_\nu\sum
\left\{\int_{\mathbb{R}^3}\left(
|\partial _{\beta _1}^{\alpha-\alpha_1}{\bf P}f|_2 |\partial _{\beta_2}^{\alpha_1}({\bf I-P})f|_\nu\right)^2
 dx\right\}^{1/2}.
\end{gather*}
Since $({\bf I-P})f$ is always part of the dissipation, we can  take the supremum in $x$ on the term with the least number of total derivatives and use the embedding $H^2(\mathbb{R}^3)\subset L^\infty(\mathbb{R}^3)$ to establish Lemma \ref{nonlin2} for this term.  And the last two terms can be estimated exactly the same as this last one.
\end{proof}

We will now introduce a bit more notation.
For $0\le m\le N$, a reduced order instant energy functional satisfies 
\begin{gather*}
\frac{1}{C}\mathcal{E}_{m}(t) \le 
\sum_{|\beta|\leq m,|\alpha|+|\beta|\leq N}
||\partial_{\beta}^{\alpha}f(t)||^{2} 
+\sum_{|\alpha |\le N}
||[\partial ^\alpha E(t),\partial ^\alpha B(t)]||^2
\le C \mathcal{E}_{m}(t). 
\end{gather*}
Similarly the reduced order dissipation rate is given by 
\begin{gather*}
\mathcal{D}_{m}(t) \equiv 
\|E(t)\|^2+\sum_{0<|\alpha |\le N}||\partial^\alpha f(t)||^2+
\sum_{|\beta|\le m, |\alpha |+|\beta |\le N}||\partial _\beta ^\alpha \{{\bf I-P}\} f(t)||^2.
\end{gather*}
Note that, $\mathcal{E}_{N}(t)=\mathcal{E}(t) $ 
and $\mathcal{D}_{N}(t)=\mathcal{D}(t)$.  
With this notation in hand, we can now prove the main global existence theorem of this last section by induction. \\

\noindent {\it Proof of Theorem \ref{energyEST}}.  We will show that for any integer $m$ with  $0\le m\le N$ we have
\begin{gather}
\frac{d}{dt} \mathcal{E}_{m}(t)+\mathcal{D}_{m}(t)\le \sqrt{{\mathcal E}(t)}\mathcal{D}(t).
\label{diffEQ}
\end{gather}
Theorem \ref{energyEST} is a special case of \eqref{diffEQ} with $m=N$.  We will prove \eqref{diffEQ}
by an induction over $m$, the order of the $v-$derivatives.
We may use different equivalent instant energy functionals on different lines without mention.

For $m=0$ and $|\alpha|\le N$, by taking the pure $\partial^\alpha$ derivatives of (\ref{vb}) we obtain 
\begin{gather}
\{\partial _t+v\cdot \nabla _x+q(E+v\times B)\cdot \nabla _v\}\partial
^\alpha f-\{\partial ^\alpha E\cdot v\}\sqrt{\mu }q_1+L\{\partial ^\alpha f\} 
\nonumber
\\
=
-
\sum_{\alpha _1\neq 0}C^\alpha_{\alpha _1}q(\partial ^{\alpha_1}E
+
v\times \partial ^{\alpha _1}B)\cdot \nabla _v\partial ^{\alpha -\alpha_1}f
 \label{purex} 
\\
+
\sum_{\alpha _1\le \alpha}C^\alpha_{\alpha _1}
\left\{\frac{q}{2}(\partial ^{\alpha _1}E\cdot v)\partial ^{\alpha -\alpha _1}f+\Gamma (\partial ^{\alpha _1}f_{,}\partial
^{\alpha -\alpha _1}f)\right\}.  \nonumber
\end{gather}
We will take the $L^2(\mathbb{R}^3_x\times \mathbb{R}^3_v)$ inner product of this with $\partial^\alpha f$.
The first term is clearly
$
\frac{1}{2}\frac d{dt}\|\partial ^\alpha f(t)\|^2.
$
From the Maxwell system \eqref{maxwell} we have
$$
-\left( \partial ^\alpha E\cdot v\sqrt{\mu }q_1,\partial ^\alpha f\right)
=
-\partial^\alpha E\cdot \partial^\alpha \mathcal{J}
=
\frac{1}{2}\frac d{dt}\left(||\partial ^\alpha E(t)||^2+||\partial ^\alpha B(t)||^2\right).
$$
Next up, since $\nu(v)\ge c(1+|v|)$, we use two applications of Cauchy-Schwartz to get
\begin{gather*}
\int_{\mathbb{R}^3\times \mathbb{R}^3}
\left| q(\partial ^{\alpha_1}E
+
v\times \partial ^{\alpha _1}B)\cdot \nabla _v\partial ^{\alpha -\alpha_1}f\right| \left| \partial ^{\alpha}f\right| dx dv
\\
\le
\int_{\mathbb{R}^3}
\left(\left| \partial ^{\alpha_1}E\right|
+
 \left|\partial ^{\alpha _1}B\right| \right)
\left\{\int_{\mathbb{R}^3} \nu(v)\left| \nabla _v\partial ^{\alpha -\alpha_1}f\right|^2 dv\right\}^{1/2} 
\left\{\int_{\mathbb{R}^3} \nu(v)\left| \partial ^{\alpha}f\right|^2 dv\right\}^{1/2}  dx 
\\
\le
\| \partial ^{\alpha}f\|_\nu
\left\{\int_{\mathbb{R}^3}
\left(\left| \partial ^{\alpha_1}E\right|
+
 \left|\partial ^{\alpha _1}B\right| \right)^2
\left\{\int_{\mathbb{R}^3} \nu(v)\left| \nabla _v\partial ^{\alpha -\alpha_1}f\right|^2 dv\right\}  dx\right\}^{1/2}.
\end{gather*}
If $|\alpha|=0$ this term is non-existent.  Therefore, for this term,
$0<|\alpha|\le N$ and $|\alpha_1|>0$ so that either $|\alpha_1|\le N/2$ or $|\alpha-\alpha_1|+1\le N/2$.  Without loss of generality say $|\alpha_1|\le N/2$; we use 
the embedding $H^2({\mathbb R}^3)\subset L^\infty ({\mathbb R}^3)$ to obtain
\begin{gather*}
\sup_{x\in\mathbb{R}^3}
\left(\left| \partial ^{\alpha_1}E\right|+\left|\partial ^{\alpha _1}B\right| \right)
\le 
C\sum_{|\bar{\alpha}|\le 2}\left(||\partial^{\bar{\alpha}} \partial ^{\alpha _1} E(t)||
+||\partial^{\bar{\alpha}} \partial^{\alpha _1} B(t)||\right)
\le \sqrt{{\mathcal E}(t)}. 
\end{gather*}
We can put the last two estimates together and use $|\alpha|>0$ to obtain
\begin{gather*}
\int_{\mathbb{R}^3\times \mathbb{R}^3}
\left| q(\partial ^{\alpha_1}E
+
v\times \partial ^{\alpha _1}B)\cdot \nabla _v\partial ^{\alpha -\alpha_1}f\right| \left| \partial ^{\alpha}f\right| dx dv
\\
\le
\| \partial ^{\alpha}f\|_\nu \sqrt{{\mathcal E}(t)}
 \| \nabla _v\partial ^{\alpha -\alpha_1}f\|_\nu
 \le
\sqrt{{\mathcal E}(t)}\mathcal{D}(t).
\end{gather*}
It remains to estimate the last two nonlinear terms on the right in \eqref{purex}.

First we estimate  $\frac{q}{2}(\partial ^{\alpha _1}E\cdot v)\partial ^{\alpha -\alpha _1}f$.  
Consider $\alpha=0$ which implies $\alpha_1=0$.  
We bound for $\| E\|$ by the dissipation to estimate this case.
We split $f$ as in \eqref{microHYDROd} so that the exponential velocity decay of the hydrodynamic part  controls the extra $|v|$ growth:
$$
\left| \left(\frac{q}{2}\{E\cdot v\}f, f\right) \right|
\le
C\int_{\mathbb{R}^3\times \mathbb{R}^3} |E| |v|\left( | {\bf P}f|^2+|({\bf I-P})f|^2 \right)dx dv.
$$
By the Sobolev embedding 
$
H^{2}(\mathbb{R}^3)\subset  L^{\infty}(\mathbb{R}^3)
$
we have
\begin{gather*}
\int_{\mathbb{R}^3\times \mathbb{R}^3} |E(t,x)| |v| |({\bf I-P})f|^2 dx dv
\le
C\left(\sup_{x\in\mathbb{R}^3}|E(t,x)|\right)\|({\bf I-P})f\|^2_\nu
\\
\le
C\left(\sum_{|\alpha|\le 2}\|\partial^\alpha E(t)\|\right)\|({\bf I-P})f\|^2_\nu
\le \sqrt{{\mathcal E}(t)}\mathcal{D}(t).
\end{gather*}
For the  term with ${\bf P}f$, we obtain
\begin{gather*}
\int_{\mathbb{R}^3\times \mathbb{R}^3} |E| |v||{\bf P}f|^2dx dv
\le 
C \left(\sup_{x\in\mathbb{R}^3} |{\bf P}f(x)|_2\right)
\int_{\mathbb{R}^3} |E| |{\bf P}f|_2dx.
\end{gather*}
By the Cauchy-Schwartz inequality  
$$
\int_{\mathbb{R}^3} |E| |{\bf P}f|_2dx
\le
C\|E(t)\|\|{\bf P} f\|
\le
\mathcal{D}^{1/2}(t) {\mathcal E}^{1/2}(t).
$$
We just controlled $\|E(t)\|$ by the dissipation, $\mathcal{D}^{1/2}(t)$.  For the other term,
from \eqref{p0}, the embedding $W^{1,6}({\mathbb R}^3)\subset L^\infty ({\mathbb R}^3)$ and the Sobolev inequality we obtain
\begin{equation}
\begin{split}
\left(\sup_{x\in\mathbb{R}^3} |{\bf P}f(x)|_2\right)
 &\le
C \|[a,b,c] \|_{W^{1,6}({\mathbb R}^3)}
 \le
C\sum_{|\alpha|\le 1} \| \nabla_x \partial^\alpha f\|
 \le C   \mathcal{D}^{1/2}(t).
\end{split}
\label{hydroSOBOLEV}
\end{equation}
Adding together the last few estimates we obtain 
\begin{equation*}
\begin{split}
\left| \left(\frac{q}{2}\{E\cdot v\}f, f\right) \right|
&\le
C\int_{\mathbb{R}^3\times \mathbb{R}^3} |E| |v|\left( | {\bf P}f|^2+|({\bf I-P})f|^2 \right)dx dv
\le
\sqrt{{\mathcal E}(t)}\mathcal{D}(t).
\end{split}
\end{equation*}
This completes the estimate when $\alpha=0$.

Alternatively if $|\alpha|>0$, we  use $\nu(v)\ge c(1+|v|)$ and Cauchy-Schwartz twice as
\begin{gather*}
\left| \left(\frac{q}{2}(\partial ^{\alpha _1}E\cdot v)\partial ^{\alpha -\alpha _1}f, \partial ^{\alpha}f\right)\right|
\le
C\int_{\mathbb{R}^3\times \mathbb{R}^3}
\left| \partial ^{\alpha_1}E\right|
\nu(v)\left| \partial ^{\alpha -\alpha_1}f\right|
\left| \partial ^{\alpha}f\right| dv  dx 
\\
\le
C\int_{\mathbb{R}^3}
\left| \partial ^{\alpha_1}E\right|
\left\{\int_{\mathbb{R}^3} \nu(v)\left| \partial ^{\alpha -\alpha_1}f\right|^2 dv\right\}^{1/2} 
\left\{\int_{\mathbb{R}^3} \nu(v)\left| \partial ^{\alpha}f\right|^2 dv\right\}^{1/2}  dx 
\\
\le
C\| \partial ^{\alpha}f\|_\nu\left\{\int_{\mathbb{R}^3}
\left| \partial ^{\alpha_1}E\right|^2
\left\{\int_{\mathbb{R}^3} \nu(v)\left| \partial ^{\alpha -\alpha_1}f\right|^2 dv\right\}
 dx\right\}^{1/2}.
\end{gather*}
Here $\| \partial ^{\alpha}f\|_\nu$ is controlled by the dissipation since $|\alpha|>0$.
But we could also have $|\alpha|=|\alpha_1|=N$ and then $\partial^\alpha E$ would not be controlled by the dissipation.    Therefore this estimate requires more care.

If $|\alpha-\alpha_1|\le N/2$, then to ensure that we have at least one derivative on $f$ we 
 use the embedding $W^{1,6}(\mathbb{R}^3)\subset L^\infty(\mathbb{R}^3)$ followed by Sobolev's inequality to obtain
\begin{gather*}
 \sup_{x\in\mathbb{R}^3}  \int_{\mathbb{R}^3} \nu(v)|\partial ^{\alpha -\alpha_1}  f(x,v)|^2dv
 \le 
 C\int_{\mathbb{R}^3} \nu(v) \| \partial ^{\alpha -\alpha_1}f(v)\|^2_{W^{1,6}(\mathbb{R}^3_x)} dv
 \\
 \le
 C\sum_{|\bar{\alpha}|\le 1}
 \|\nabla_x \partial ^{\alpha -\alpha_1}\partial^{\bar{\alpha}}f\|_\nu^2
 \le 
 C \mathcal{D}(t).
 \end{gather*}
 Alternatively, if $|\alpha-\alpha_1|> N/2$ then $|\alpha_1|\le N/2$ and 
 \begin{gather}
 \sup_{x\in\mathbb{R}^3}  |\partial ^{\alpha_1}  E(t,x)|^2
 \le 
 C \| \partial ^{\alpha_1}E\|^2_{W^{1,6}(\mathbb{R}^3_x)} 
 \le
 C\sum_{|\bar{\alpha}|\le 1}
 \|\nabla_x \partial ^{\alpha_1}\partial^{\bar{\alpha}}E\|^2
 \le 
 C \mathcal{E}(t).
 \label{Eembed}
 \end{gather}
In either case the last two estimates establish
\begin{gather*}
\int_{\mathbb{R}^3}
\left| \partial ^{\alpha_1}E\right|^2
\left\{\int_{\mathbb{R}^3} \nu(v)\left| \partial ^{\alpha -\alpha_1}f\right|^2 dv\right\}
 dx
  \le
\mathcal{E}(t) \mathcal{D}(t).
\end{gather*}
Combining these last few estimates, since $|\alpha|>0$, we obtain
\begin{gather*}
\left| \left(\frac{q}{2}(\partial ^{\alpha _1}E\cdot v)\partial ^{\alpha -\alpha _1}f, \partial ^{\alpha}f\right)\right|
\le
\| \partial ^{\alpha}f\|_\nu
\sqrt{{\mathcal E}(t)} \mathcal{D}^{1/2}(t)
\le
\sqrt{{\mathcal E}(t)} \mathcal{D}(t).
\end{gather*}
It remains to estimate the nonlinear collision operator in \eqref{purex}.

But from Lemma \ref{nonlin2} and the basic invariant property \eqref{invariantP} we have
$$
\left| \left( \Gamma (\partial^{\alpha -\alpha _{1}}f,\partial ^{\alpha_{1}}f),\partial ^{\alpha }f\right) \right| 
\leq 
\sqrt{{\mathcal E}(t)} \mathcal{D}(t).
$$
Adding up all the estimates for the terms in \eqref{purex} and summing over $|\alpha|\le N$ yields
\begin{equation}
\frac{1}{2}\frac{d}{dt}\left\{ \sum_{|\alpha|\le N}
\|\partial ^{\alpha}f\|^{2}+\|\partial ^{\alpha}[E,B]\|^{2}\right\} 
+
\sum_{|\alpha|\le N}\left( L\partial^\alpha f(t), \partial^\alpha f(t)\right)
\leq
\sqrt{{\mathcal E}(t)} \mathcal{D}(t).
\label{step0}
\end{equation}
By Theorem \ref{positive} and  Lemma \ref{field} for some small $\delta_{0}^\prime>0$ we see that
\begin{gather*}
\frac{1}{2}\frac{d}{dt}\left\{ \sum_{|\alpha|\le N}
\|\partial ^{\alpha}f\|^{2}+\|\partial ^{\alpha}[E,B]\|^{2}-2C_{0}G(t)\right\} 
+
\delta_{0}^\prime \mathcal{D}_0(t)
\leq
\sqrt{{\mathcal E}(t)} \mathcal{D}(t).
\end{gather*}
We add \eqref{step0} multiplied by a large constant, $C_0^\prime>0$, to the last line to obtain
\begin{gather}
\frac{d}{dt}{\mathcal E}_0(t)
+ \mathcal{D}_0(t)
\leq
\sqrt{{\mathcal E}(t)} \mathcal{D}(t).
\label{lzero}
\end{gather}
We have used $L\ge 0$ from Lemma \ref{llinear}.  Above, for $C_0^\prime$ large enough and from the definition of $G(t)$ in Theorem \ref{positive}, the following is a reduced instant energy functional:
$$
\mathcal{E}_0(t)
=
\frac{(C_0^\prime \delta_{0}^\prime +1)}{2 \delta_{0}^\prime}
\left\{\sum_{|\alpha|\le N}
\|\partial ^{\alpha}f\|^{2}+\|[\partial ^{\alpha} E, \partial ^{\alpha} B]\|^{2}\right\}-\frac{C_{0}}{\delta_{0}^\prime}G(t).
$$
We conclude \eqref{diffEQ} for $\beta=0$.

Say Theorem \ref{positive} holds for $|\beta |=m>0$. For $|\beta |=m+1$,  take $\partial _\beta ^\alpha$ of 
(\ref{vb}) to obtain: 
\begin{gather}
\{\partial _t+v\cdot \nabla _x+q(E+v\times B)\cdot \nabla _v\}\partial_\beta ^\alpha ({\bf I-P})f
-\partial ^\alpha E\cdot \partial _\beta \{v\sqrt{\mu }\}q_1
 \nonumber 
\\
+\partial _\beta \{L\partial ^\alpha ({\bf I-P})f\} 
+\sum_{\beta _1\neq 0}C^\beta_{\beta _1}
\partial _{\beta _1}v\cdot \nabla_x\partial _{\beta -\beta _1}^\alpha f  
\nonumber 
\\
=
\sum C^\alpha_{\alpha _1}C^\beta_{\beta _1}\frac{q}{2}\{\partial ^{\alpha_1}
E\cdot \partial _{\beta _1}v\}\partial _{\beta -\beta _1}^{\alpha -\alpha_1}f
+
\partial _\beta^\alpha 
\Gamma (f, f)
\label{v} 
\\
-
\sum_{|\alpha _1|\neq 0}C^\alpha_{\alpha _1}q\partial ^{\alpha_1}
E\cdot \nabla _v\partial _\beta ^{\alpha -\alpha _1}f  
-
\sum_{|\alpha _1|+|\beta _1|\neq 0}C^\alpha_{\alpha _1}C^\beta_{\beta _1}
q\partial _{\beta _1}v\times \partial ^{\alpha _1}B\cdot \nabla_v
\partial_{\beta -\beta _1}^{\alpha -\alpha _1}f
\nonumber
\\
\nonumber
-
\{\partial _t+v\cdot \nabla _x+q(E+v\times B)\cdot \nabla _v\}\partial_\beta ^\alpha {\bf P}f.
\end{gather}
Notice that we have split the equation in terms of microscopic and hydrodynamic parts in a different way from \eqref{macro}.  We do this to estimate the linear part of the collision operator properly in terms of our weak whole space dissipation.  Standard estimates for $(\partial_\beta L \partial^\alpha f, \partial_\beta^\alpha f)$, e.g. Lemma \ref{kcompact}, are too strong in this case.

We take the inner product of (\ref{v}) with $\partial _\beta ^\alpha ({\bf I-P})f$ over 
${\mathbb R}^3_x\times {\mathbb R}^3_v$.  We estimate each term  from left to right.  The first inner product on the left gives 
\[
\frac 12\frac d{dt}||\partial _\beta ^\alpha ({\bf I-P}) f(t)||^2.
\]
Since $|\beta |=m+1>0,$ $|\alpha |\le N-m-1$.  Then by Lemma \ref{field} the
second inner product on the left is bounded by
\begin{gather*}
-\left| \left(  \partial ^\alpha E\cdot \partial _\beta \{v\sqrt{\mu }%
\}q_1,\partial _\beta ^\alpha ({\bf I-P})f \right)  \right|
\ge 
-C \|\partial ^\alpha E\|  \|\partial_\beta^\alpha ({\bf I-P})f \|
\\
\ge 
-C_\eta \|\partial ^\alpha E\|^2
-\eta \|\partial_\beta^\alpha ({\bf I-P})f \|^2
\\
\ge
-C_\eta \mathcal{D}_0(t)-\eta \|\partial_\beta^\alpha ({\bf I-P})f \|^2.
\end{gather*}
For the linear operator,  since $|\partial_{\beta _1}\nu (v)|$ is bounded for $\beta _1\ne 0$ (Lemma \ref{kcompact}) we have 
\begin{gather*}
\left( \partial _\beta \{\nu \partial ^\alpha ({\bf I-P})f\},\partial _\beta ^\alpha ({\bf I-P})f\right)
=
||\partial _\beta ^\alpha({\bf I-P}) f||_\nu ^2
\\
+
\sum_{\beta_1\ne 0} C^\beta_{\beta_1}\langle  \partial _{\beta_1} \nu \partial_{\beta-\beta_1}^\alpha ({\bf I-P})f,\partial _\beta ^\alpha ({\bf I-P})f\rangle 
\\
\ge ||\partial _\beta ^\alpha ({\bf I-P}) f||_\nu ^2
-
C||\partial _\beta ^\alpha ({\bf I-P}) f||\sum_{\beta_1\ne 0} 
||\partial _{\beta-\beta _1}^\alpha ({\bf I-P}) f||.
\end{gather*}
Then together with Lemma \ref{kcompact}, since $L=\nu-K$, we deduce that for any $\eta >0$ there
is a constant $C_\eta >0$ such that the third term on the left side of \eqref{v} is bounded from
below as 
\begin{gather*}
\left(  \partial _\beta \{L\partial^\alpha ({\bf I-P})f\},\partial _\beta ^\alpha
({\bf I-P})f  \right)
 \ge 
 ||\partial _\beta ^\alpha ({\bf I-P}) f||_\nu ^2
-
\eta \sum_{|\beta^\prime|=|\beta |}\|\partial _{\beta^\prime}^\alpha ({\bf I-P})f\|_\nu ^2
\\
-C_\eta
\sum_{|\bar{\beta}|<|\beta |}||\partial_{\bar{\beta}}^\alpha ({\bf I-P}) f||^2.
\end{gather*}
We can further choose $C_\eta >0$ such that the inner product of the
last term on left side of (\ref{v}) is  bounded by 
\begin{gather*}
\eta ||\partial _\beta ^\alpha ({\bf I-P}) f(t)||^2+C_\eta \sum_{|\beta _1|=1}||\nabla
_x\partial _{\beta -\beta _1}^\alpha f||^2.
\end{gather*}
We split the last term above as in \eqref{microHYDROd} to obtain
\begin{gather*}
\le
\eta ||\partial _\beta ^\alpha ({\bf I-P}) f(t)||^2+C_\eta \sum_{|\beta _1|=1}||\nabla
_x\partial _{\beta -\beta _1}^\alpha {\bf P} f||^2
+C_\eta \sum_{|\beta _1|=1}||\nabla
_x\partial _{\beta -\beta _1}^\alpha ({\bf I-P}) f||^2
\\
\le
\eta ||\partial _\beta ^\alpha ({\bf I-P}) f(t)||^2+C\mathcal{D}_{0}(t)
+C_\eta \sum_{|\beta _1|=1}||\nabla
_x\partial _{\beta -\beta _1}^\alpha ({\bf I-P}) f||^2.
\end{gather*}
This splitting was used to get rid of the velocity derivatives of the hydrodynamic part, which are not in the dissipation.

Next we estimate all the terms on the right side of \eqref{v}.  Let us first consider the first term on the right, using $(1+|v|)\le c \nu(v)$ we have
\begin{gather}
\label{111}
\left| \left( \frac{q}{2}\{\partial ^{\alpha_1}
E\cdot \partial _{\beta _1}v\}\partial _{\beta -\beta _1}^{\alpha -\alpha_1}f,
\partial _{\beta}^{\alpha}({\bf I-P})f \right) \right|
\\
\le
C\int_{\mathbb{R}^3\times \mathbb{R}^3}
\left|\partial ^{\alpha_1}E\right| \nu(v)
\left| \partial _{\beta -\beta _1}^{\alpha -\alpha_1}f \right|
\left| \partial _{\beta}^{\alpha}({\bf I-P})f \right| dx dv.
\nonumber
\end{gather}
We use Cauchy-Schwartz twice to see that the last line is
\begin{gather*}
\le
C\int_{\mathbb{R}^3}
\left|\partial ^{\alpha_1}E\right| \left\{
\int_{\mathbb{R}^3} \nu(v)
\left| \partial _{\beta -\beta _1}^{\alpha -\alpha_1}f \right|^2 dv\right\}^{1/2}
\left\{\int_{\mathbb{R}^3}\nu(v)\left| \partial _{\beta}^{\alpha}({\bf I-P}) f \right|^2 dv\right\}^{1/2} dx 
\\
\le
C\|\partial _{\beta}^{\alpha} ({\bf I-P})f \|_\nu  \left\{\int_{\mathbb{R}^3}
\left|\partial ^{\alpha_1}E\right|^2 \left\{
\int_{\mathbb{R}^3} \nu(v)
\left| \partial _{\beta -\beta _1}^{\alpha -\alpha_1}f \right|^2 dv\right\}dx \right\}^{1/2}.
\end{gather*}
We split $f$ into it's microscopic and hydrodynamic parts  \eqref{microHYDROd} so that we can remove the velocity derivatives from the hydrodynamic part:
\begin{gather*}
\le
C\|\partial _{\beta}^{\alpha} ({\bf I-P})f \|_\nu  \left\{\int_{\mathbb{R}^3}
\left|\partial ^{\alpha_1}E\right|^2 \left\{
\int_{\mathbb{R}^3} \nu(v)
\left| \partial _{\beta -\beta _1}^{\alpha -\alpha_1}{\bf P}f \right|^2 dv\right\}dx \right\}^{1/2}
\\
+
C\|\partial _{\beta}^{\alpha} ({\bf I-P})f \|_\nu  \left\{\int_{\mathbb{R}^3}
\left|\partial ^{\alpha_1}E\right|^2 \left\{
\int_{\mathbb{R}^3} \nu(v)
\left| \partial _{\beta -\beta _1}^{\alpha -\alpha_1}({\bf I-P}) f \right|^2 dv\right\}dx \right\}^{1/2}.
\end{gather*}
Since the total number of derivatives is at most $N$, either $|\alpha_1|\le N/2$ or $|\alpha-\alpha_1|+|\beta-\beta_1|\le N/2$.
So we take the supremum in $x$ of the term has the least total derivatives, for the microscopic term, using the embedding $H^2(\mathbb{R}^3) \subset L^\infty(\mathbb{R}^3)$ to observe
\begin{gather*}
\int_{\mathbb{R}^3}
\left|\partial ^{\alpha_1}E\right|^2 \left\{
\int_{\mathbb{R}^3} \nu(v)
\left| \partial _{\beta -\beta _1}^{\alpha -\alpha_1}({\bf I-P}) f \right|^2 dv\right\}dx 
\le 
{\mathcal E}(t) \mathcal{D}(t).
\end{gather*}
On the other hand, for the hydrodynamic variables we have
\begin{gather*}
\int_{\mathbb{R}^3}
\left|\partial ^{\alpha_1}E\right|^2 \left\{
\int_{\mathbb{R}^3} \nu(v)
\left| \partial _{\beta -\beta _1}^{\alpha -\alpha_1}{\bf P}f \right|^2 dv\right\}dx 
\le
C
\int_{\mathbb{R}^3}
\left|\partial ^{\alpha_1}E\right|^2 
|\partial^{\alpha -\alpha_1}{\bf P}f |_2^2 dx.
\end{gather*}
If $|\alpha-\alpha_1|\le N/2$, then as in \eqref{hydroSOBOLEV} we see that 
$$
\int_{\mathbb{R}^3}
\left|\partial ^{\alpha_1}E\right|^2 
|\partial^{\alpha -\alpha_1}{\bf P}f |_2^2 dx
\le
C\|\partial ^{\alpha_1}E\|^2 
\mathcal{D}(t)
\le
C
\mathcal{E}(t)
\mathcal{D}(t).
$$
If alternatively $|\alpha-\alpha_1|> N/2$ and $|\alpha_1|\le N/2$ then we use \eqref{Eembed} to get
$$
\int_{\mathbb{R}^3}
\left|\partial ^{\alpha_1}E\right|^2 
|\partial^{\alpha -\alpha_1}{\bf P}f |_2^2 dx
\le
C
\mathcal{E}(t)
\|\partial^{\alpha -\alpha_1}f \|_\nu^2
\le
C
\mathcal{E}(t)
\mathcal{D}(t).
$$
In any case, we conclude
\begin{gather*}
\left| \left( \frac{q}{2}\{\partial ^{\alpha_1}
E\cdot \partial _{\beta _1}v\}\partial _{\beta -\beta _1}^{\alpha -\alpha_1}f,
\partial _{\beta}^{\alpha}({\bf I-P})f \right) \right|
\le
\sqrt{{\mathcal E}(t)} \mathcal{D}(t).
\end{gather*}
This is the estimate for the first term on the right side of \eqref{v}.

For the nonlinear collision operator in \eqref{v},  by Lemma \ref{nonlin2} we have
$$
\left| \left(  \partial _\beta^\alpha 
\Gamma (f, f),
\partial _{\beta}^{\alpha} ({\bf I-P})f\right) \right|   
\le 
\sqrt{\mathcal{E}(t)}\mathcal{D}(t).
$$
For the next two  terms on the right side of \eqref{v} we have 
\begin{gather*}
\left| \left(
q\partial ^{\alpha_1}
E\cdot \nabla _v\partial _\beta ^{\alpha -\alpha _1}f
+
 q\partial _{\beta _1}v\times \partial ^{\alpha _1}B\cdot \nabla_v
\partial_{\beta -\beta _1}^{\alpha -\alpha _1}f, \partial_{\beta}^{\alpha}({\bf I-P})f \right) \right| 
\\
\le
C
\int_{\mathbb{R}^3\times \mathbb{R}^3} \nu(v)\left(  |\partial^{\alpha_1}E| 
\left|\nabla_v  \partial_{\beta}^{\alpha -\alpha _1}f \right|
+ 
\left|\partial^{\alpha _1}B\right|
\left| \nabla_v \partial_{\beta -\beta _1}^{\alpha -\alpha _1}f \right| 
\right)
\left|  \partial_{\beta }^{\alpha }({\bf I-P})f \right|  dv dx.
\end{gather*}
Notice that these are in the form \eqref{111}.  They are therefore estimated the same way.

Furthermore the  same procedure (as in the estimate for \eqref{111}) yields
$$
\left| \left(q(E+v\times B)\cdot \nabla _v\}\partial_\beta ^\alpha {\bf P}f, 
\partial_\beta ^\alpha ({\bf I-P})f \right) \right|
\le \sqrt{\mathcal{E}(t)}\mathcal{D}(t).
$$
It remains to estimate 
\begin{gather*}
\left| \left(
\{\partial _t+v\cdot \nabla _x\}\partial_\beta ^\alpha {\bf P}f,
\partial_\beta ^\alpha ({\bf I-P})f \right) \right|.
\end{gather*}
Since $|\beta|>0$, we have $|\alpha| \le N-1$ and the above is
\begin{gather*}
\le
C\int_{\mathbb{R}^3\times\mathbb{R}^3} 
\nu(v)\left( \left|\partial _t\partial_\beta ^\alpha {\bf P}f\right|+\left|\nabla _x \partial_\beta ^\alpha {\bf P}f\right|\right)
\left|\partial_\beta ^\alpha ({\bf I-P})f \right| dx dv
\\
\le
C\|\partial_\beta ^\alpha ({\bf I-P})f \|_\nu 
\left( \|\partial _t\partial_\beta ^\alpha {\bf P}f\|_\nu +\|\nabla _x \partial_\beta ^\alpha {\bf P}f\|_\nu \right)
\\
\le
C\|\partial_\beta ^\alpha ({\bf I-P})f \|_\nu 
\left( \|\partial _t\partial^\alpha {\bf P}f\| +\|\nabla _x \partial^\alpha {\bf P}f\| \right)
\\
\le
\eta \|\partial_\beta ^\alpha ({\bf I-P})f \|_\nu^2 
+C_\eta \sum_{|\bar{\alpha}|= 1}\|\partial^{\bar{\alpha}}\partial^\alpha f\|^2.
\end{gather*}
This completes  all the estimates for the terms in (\ref{v}).

By collecting all the estimates for the terms in \eqref{v} and summing over 
$|\beta |=m+1$ and $|\alpha |+|\beta |\le N$ we obtain
\begin{gather*}
\sum_{|\beta |=m+1, |\alpha |+|\beta |\le N} 
\left\{ \frac{1}{2}\frac{d}{dt}\|\partial _\beta ^\alpha ({\bf I-P}) f(t)\|^2
+||\partial _\beta ^\alpha ({\bf I-P}) f(t)||_\nu ^2\right\}  -C_\eta \mathcal{D}_0(t)
\\
\le 
C \eta
\sum_{|\beta |=m+1,
|\alpha |+|\beta |\le N}
||\partial _\beta ^\alpha ({\bf I-P})f(t)||_\nu ^2
\\
+
C_\eta
\sum_{|\beta |\le m, |\alpha |+|\beta |\le N}
||\partial _\beta ^\alpha ({\bf I-P}) f(t)||^2
+
\sqrt{{\mathcal E}(t)}{\mathcal D}(t).
\end{gather*}
Here $C$ is a large constant which does not depend on $\eta$.  Choosing 
$
\eta=\frac{1}{2C}
$
we have
\begin{gather}
\sum_{|\beta |=m+1, |\alpha |+|\beta |\le N} 
\frac{1}{2}\left\{ \frac{d}{dt}\|\partial _\beta ^\alpha ({\bf I-P}) f(t)\|^2
+||\partial _\beta ^\alpha ({\bf I-P}) f(t)||_\nu ^2\right\}  
\nonumber
\\
\le 
C
\mathcal{D}_m(t)
+
\sqrt{{\mathcal E}(t)}{\mathcal D}(t).
\label{v1}
\end{gather}
Now add the inequality from \eqref{diffEQ} for $|\beta|=m$ multiplied by a suitably large constant, $C_m$, to \eqref{v1} to obtain
\begin{gather*}
\frac{d}{dt}\left\{ C_m\mathcal{E}_m(t)
+\frac{1}{2}\sum_{|\beta |=m+1, |\alpha |+|\beta |\le N} \|\partial _\beta ^\alpha ({\bf I-P}) f(t)\|^2\right\}
+
\left(C_m-C\right)\mathcal{D}_m(t) 
\\
+
\frac{1}{2}\sum_{|\beta |=m+1, |\alpha |+|\beta |\le N} 
\|\partial _\beta ^\alpha ({\bf I-P}) f(t)\|_\nu ^2
\le 
\sqrt{{\mathcal E}(t)}{\mathcal D}(t).
\end{gather*}
Here $C_m$ is chosen so that $C_m-C>0$.
We further define
\begin{gather*}
\mathcal{E}_{m+1}(t)
=
C_m\mathcal{E}_m(t)
+\frac{1}{2}\sum_{|\beta |=m+1, |\alpha |+|\beta |\le N} \|\partial _\beta ^\alpha ({\bf I-P}) f(t)\|^2.
\end{gather*}
Since $\| \partial_\beta {\bf P}f \| \le C\|  {\bf P}f \|$, $\mathcal{E}_{m+1}(t)$
is an instant energy functional.
\qed

\bigskip

\noindent {\bf Acknowledgements}.  The author is very grateful to  Yan Guo for bringing his attention to this problem.
He thanks the referees for several helpful comments concerning the presentation of this paper, and for pointing out the reference \cite{MR1997264}.
He also thanks Tong Yang for kindly sending \cite{YZ2005} in December 2005.

\end{document}